\documentclass[preprint,12pt]{elsarticle}

\usepackage{amsthm}
\usepackage{longtable}
\usepackage{bbding}
\usepackage{amsmath}
\usepackage{amssymb}
\usepackage{multirow}
\usepackage{booktabs}
\usepackage{colortbl}
\usepackage{fancyhdr}
\usepackage{indentfirst}
\usepackage{amsmath}
\usepackage{makecell}
\usepackage{multirow}
\usepackage{array, threeparttable}
\usepackage{cite}
\usepackage{url}
\usepackage{breakurl}
\usepackage{footnote}
\usepackage{threeparttable}
\usepackage{hyperref}
\usepackage{cleveref}
\usepackage{setspace}
\usepackage{lineno,hyperref}
\usepackage{ifthen,url,graphicx,color,array,cases}
\usepackage{geometry}
\usepackage{amssymb}
\newtheorem{theorem}{Theorem}
\newtheorem{proposition}{Proposition}
\newtheorem{lemma}{Lemma}
\newtheorem{remark}{Remark}
\usepackage{caption}
\usepackage{algorithm, algorithmic}
\usepackage{appendix}

\allowdisplaybreaks[4]

\journal{Transportation Research Part E: Logistics and Transportation Review.~~}

\usepackage{setspace}

\begin{document}


\begin{frontmatter}
\title{A Benders and column generation method to the integrated airline schedule and aircraft recovery with gate reassignment}

\author[1]{Jianlin Jiang}\ead{jiangjianlin@nuaa.edu.cn}
\author[2]{Jiajin Lin}\ead{jiajin\_lin@outlook.com}
\author[1]{Yan Gu\corref{corauthor}}\ead{guyanmath@nuaa.edu.cn}
\author[3]{Yibing Lv}\ead{Yibinglv@yangtzeu.edu.cn}
\author[4]{Cheng-Lung Wu}\ead{c.l.wu@unsw.edu.au}

\address[1]{School of Mathematics,
            Nanjing University of Aeronautics and Astronautics,\\
            Nanjing 210016,
            China}
\address[2]{School of Information Engineering,
            Taizhou Vocational College of Science and Technology,\\
            Taizhou 318020,
            China}
\address[3]{School of Information and Mathematics, Yangtze University,\\
            Jingzhou 434023,
            China}
\address[4]{School of Aviation, University of New South Wales,\\
            Sydney, Kensington, NSW 2052,
            Australia}
\cortext[corauthor]{Corresponding author.}


\begin{abstract}
Disruptions are inevitable during airline operations, and disruptions cost airlines and the traveling public. Disruption recovery decisions are often made in a sequential process that includes flight rescheduling, aircraft rerouting, crew reassignment, passenger re-accommodation, and gate reassignment for airline operations. Such a sequential recovery approach alleviates the complexity of the whole recovery process but often causes further disruptions to airport gate assignments and high recovery costs for airlines and passengers. This paper presents a novel airline disruption recovery approach by integrating the schedule and aircraft recovery with gate reassignment. We propose a Benders and column generation (BCG) method to solve the integrated problem. By exploiting the favorable structure arising from the Benders decomposition framework, we provide two acceleration techniques for Benders subproblems, a separation technique and an effective infeasibility certificate, to enhance the efficiency of the BCG method. The proposed model is tested on real-world data. In our experiments, all test instances are solved by the BCG method under a five-minute threshold with optimality gaps within 5\%. Compared with the sequential recovery approach, the integrated recovery approach significantly saves recovery costs and avoids infeasible gate reassignments caused by the sequential recovery approach.
\end{abstract}


\end{frontmatter}


\section{Introduction}\label{Introduction}
Robust scheduling has been demonstrated as a key concept in the literature to proactively improve the resilience of airline schedule execution (\citeauthor{dunbar2012robust}, \citeyear{dunbar2012robust}; \citeauthor{ahmed2018robust}, \citeyear{ahmed2018robust}; \citeauthor{wu2018airlinecapacity}, \citeyear{wu2018airlinecapacity}). However, unpredictable disruptions are inevitable in daily airline operations. These disruptions can cause flight delays, cancellations, and travel inconveniences to passengers. When disruptions happen, the airline operations control center starts the {\it reactive} disruption management process by reallocating resources (such as flights, aircraft, and crews) to recover the disrupted schedules. Then, an efficient and effective recovery process can significantly lower the operating loss to an airline while maintaining its responsibility for transporting passengers.\par

In general, reactive airline disruption management is divided into several phases: schedule recovery (timetable recovery), aircraft recovery (aircraft rerouting), crew recovery (crew reassignment), and passenger recovery (passenger re-accommodation) (\citeauthor{wu2018airlinecapacity}, \citeyear{wu2018airlinecapacity}). Gate reassignment is also an important phase in disruption management for both airlines and airports. These recovery phases are often conducted sequentially in a tight time frame due to the complexity of the disruption management problem. The sequential approach simplifies the recovery process, but this approach ignores the interdependence among these recovery phases. This may result in sub-optimality for the recovery solution or infeasibility to subsequent recovery phases. In the last two decades, there has been increasing interest in integrating two or more recovery phases to improve solution quality and avoid infeasibility (see for example, \citeauthor{petersen2012optimization}, \citeyear{petersen2012optimization}; \citeauthor{maher2016solving}, \citeyear{maher2016solving}). For more details about disruption management, please refer to \citeauthor{wu2025} (\citeyear{wu2025}) for a recent literature review.\par

This paper focuses on an integrated recovery problem, including schedule recovery, aircraft recovery, and gate reassignment. The following three {\it motivations} promote us to explore the integration of these phases. First, such integration can avoid further disruptions to gate reassignment due to the uncertainty in estimating gate capacity in a sequential recovery approach. To implicitly consider the gate capacity when solving the schedule and aircraft recovery problems, airlines adjust the slot capacity forecasts (the maximum allowable number of departure/arrival flights in slots). Since forecasts are highly uncertain during disruptions, airport gate reassignments are often further disrupted and require recovery (\citeauthor{kontoyiannakis2009simulation}, \citeyear{kontoyiannakis2009simulation}; \citeauthor{yan2016tarmac}, \citeyear{yan2016tarmac}). When the slot capacity is overestimated, delayed flights from the previous slots may still occupy gates, and arriving flights during the current slot may not have available gates upon arrival, causing {\it gate blockage} (\citeauthor{castaing2016reducing}, \citeyear{castaing2016reducing}; \citeauthor{dorndorf2017reducing}, \citeyear{dorndorf2017reducing}). When the slot capacity is underestimated, the risk of gate blockage can be mitigated by the conservative slot capacity forecast, but more flights will be delayed and available gates are underutilized, leading to a high recovery cost for airlines and airports (\citeauthor{flyertalk2019}, \citeyear{flyertalk2019}). The above two issues indeed exist in airport and airline operations (\citeauthor{PYOK2024}, \citeyear{PYOK2024}; \citeauthor{WLOS2024}, \citeyear{WLOS2024}). To both avoid gate blockage and reduce flight delays, this paper explores the benefits and feasibility of integrating gate reassignment in the schedule and aircraft recovery problem to overcome the possible operational issues caused by the slot capacity forecasts in a sequential recovery approach. Computational experiments in Section \hyperref[Section5]{6} to illustrate this idea further.\par

Second, the integration of these phases can improve the overall recovery performance by capturing the network effect. In practical operations and the existing literature, given an aircraft recovery plan, gate reassignment is usually conducted at the airport level, i.e., at individual airports (\citeauthor{yan2011airport}, \citeyear{yan2011airport}; \citeauthor{zhang2017optimization}, \citeyear{zhang2017optimization}). This approach alleviates the complexity of the airline recovery process but ignores the network effect in disruption management. Specifically, gate delays and gate blockage at one airport will disrupt other flights at that airport due to the resource scheduling and subsequent flights on the associated aircraft routes at other airports. This further affects gate assignment in the airline network, which causes a ripple effect. The worst scenario is that an airline may need to readjust its previous recovery decisions later because of the network delay propagation. Based on this rationale, it is hypothetically beneficial to explore the benefits of integrating gate reassignment with the schedule and aircraft recovery to further capture the interdependence among them in an airline network.\par

Third, the proposed integration can facilitate airport collaborative decision-making (CDM) during disruption recovery, an initiative by the International Civil Aviation Organization (ICAO). In the industry, airport gates are managed by airport authorities or airlines (\citeauthor{zhang2017optimization}, \citeyear{zhang2017optimization}). In the context of gate management under the airport authority, to mitigate disruptions as soon as possible, CDM adjust airport resources such as slots and gates efficiently to achieve the best outcome for stakeholders (\citeauthor{geske2024}, \citeyear{geske2024}). As such, the integration of the schedule and aircraft recovery and the gate reassignment can realize this idea. When airlines manage their own gates at an airport (\citeauthor{tang2013airport}, \citeyear{tang2013airport}), they can enhance recovery decisions by considering gate constraints and operational requirements in schedule adjustment and aircraft routing recovery. In this case, airlines can both avoid gate blockage and save recovery costs in operations.  From the above two viewpoints, there are potential benefits in exploring an integrated approach for flight rescheduling, aircraft rerouting, and gate reassignment.\par

The main contributions of this paper are highlighted as follows:
\vspace{-0.2cm}
\begin{itemize}
\setlength{\itemsep}{2pt}
\setlength{\parsep}{0pt}
\setlength{\parskip}{0pt}

\item We develop an integrated schedule, aircraft and gate recovery model (SAGRM) to formulate the airline disruption recovery problem. An integrated approach for recovering these phases improves recovery performance and avoids the infeasibility issue when solving the gate reassignment problem in a sequential approach. This approach reduces the recovery costs in disruption management.

\item A Benders and column generation (BCG) method is proposed to solve the integrated recovery problem with a column generation (CG) algorithm to avoid enumerating all aircraft routes and gate patterns. Unlike other CG processes in the literature, we use logical expressions to generate aircraft routes and gate patterns by using a constraint programming (CP) model (see \ref{appendix CP models}).  Commercial CP solvers such as CPLEX perform well for solving the CP model. This approach alleviates the heavy burden of generating aircraft routes and gate patterns in the CG algorithm.

\item Based on the special feature of the SAGRM formulation, we utilize its variable separation characteristic to separate the Benders subproblem (BSP, a gate reassignment problem for all airports in an airline network) into a number of gate reassignment problems for each available gate type at each airport. These smaller parallel problems are solved by using parallel computing to significantly reduce the runtime of SAGRM for real-time disruption recovery needs.

\item Since not all decision variables are in the restricted BSP during the CG process, checking the feasibility of the restricted BSP is not a trivial task. We propose a generalized certificate based on \citeauthor{petersen2012optimization} (\citeyear{petersen2012optimization}) to detect the infeasibility of our separated BSP. Once the infeasibility of one separated subproblem is detected, the CG process for checking the feasibility of the BSP can be terminated immediately. It is noted that this certificate is applicable to the Benders subproblem with capacity inequality constraints, rather than only considering the linear equality constraints.\vspace{-0.2cm}
\end{itemize}

The remainder of this paper is organized as follows. Section \ref{Literature review} provides the literature review about airline disruption management. Section \hyperref[Section2]{3} develops the formulation of SAGRM to describe the integrated recovery problem. A BCG method is developed in Section \hyperref[Section3]{4}. Based on the structure of the integrated recovery problem, Section \hyperref[Section4]{5} proposes two acceleration techniques for the BCG method. Section \hyperref[Section5]{6} presents computational results to demonstrate the performance of the integrated recovery approach. The last section concludes our work and provides future work.\par

\section{Literature review}\label{Literature review}
Airline recovery is significant for the airline industry and deeply affects airlines' daily operations. This problem is usually intractable due to the complex restrictions of available resources and the real-time requirement. \citeauthor{clausen2010disruption} (\citeyear{clausen2010disruption}) summarized the previous literature on airline recovery and emphasized the similarities between the solution approaches applied to disruption management and operational planning. \citeauthor{hassan2021airline} (\citeyear{hassan2021airline}) and \citeauthor{su2021airline} (\citeyear{su2021airline}) provided a comprehensive literature survey on airline disruption management from 2009 to 2020. More recently, \citeauthor{wu2025} (\citeyear{wu2025}) further summarized the literature on airline recovery from 2021 to 2024. The authors found that airline disruption management is a growing field of research. In the following, we focus the literature review on individual recovery phases and the integrated airline recovery problem.\par

\subsection{Individual recovery phases}
\textit{Schedule recovery problem} seeks to find recovery options to repair a flight schedule after disruptions. \citeauthor{teodorovic1984optimal} (\citeyear{teodorovic1984optimal}) provided the pioneering work that repaired a disrupted schedule and obtained new aircraft routes. The proposed network flow models were solved by a heuristic branch-and-bound method. Based on a heuristic and sequential approach, \citeauthor{teodorovic1990model} (\citeyear{teodorovic1990model}, \citeyear{teodorovic1995model}) further considered flight cancellations as a recovery option and added aircraft maintenance and crew in their proposed mechanism. \citeauthor{thengvall2001multiple} (\citeyear{thengvall2001multiple}) provided three network flow models to describe a multi-fleet schedule recovery problem. Due to the special structure of the arc-based models in \citeauthor{thengvall2001multiple} (\citeyear{thengvall2001multiple}), their multi-fleet problems can be separated into several single-fleet ones. The arc-based approach can explicitly describe all practical requirements by introducing a lot of variables and constraints, but it imposes a great burden on solving the associated recovery model. \citeauthor{wang2019} (\citeyear{wang2019}) proposed a simulation-based approach to deal with schedule recovery problem. They used a push-back strategy to quickly obtain a recovery plan at the cost of the solution quality. Unlike the above work, we apply the path-based modelling approach in this paper to overcome the difficulty of finding the solution. Then, with the aid of column generation technique, we can obtain high-quality solutions while meeting the real-time requirements in disruption management.\par

\textit{Aircraft recovery problem} reassigns aircraft routes that satisfy maintenance requirements to accommodate the repaired schedule. \citeauthor{yan1996decision} (\citeyear{yan1996decision}) formulated a single-fleet aircraft recovery problem as an arc-based model on a time-space network and solved it by Lagrangian relaxation. \citeauthor{rosenberger2003rerouting} (\citeyear{rosenberger2003rerouting}) introduced cancellation cycles and single swap routes to heuristically identify disruptable aircraft that need to be considered and then developed a path-based model. \citeauthor{liang2018column} (\citeyear{liang2018column}) regarded planned maintenance activities as special flights and formulated a path-based model on a connection network with maintenance flexibility. \citeauthor{zang2024} (\citeyear{zang2024}) proposed a proactive aircraft recovery approach and implemented it with a decision-decomposition-based algorithm based on collaborative allocation of airport capacity and flow.\par

\textit{Gate reassignment problem} adjusts the pre-assignment gate schedule to accommodate the changes in the flight schedule. \citeauthor{yan2011airport} (\citeyear{yan2011airport}) extended the work of \citeauthor{yan2009airport} (\citeyear{yan2009airport}) and developed a dynamic gate reassignment approach based on a dynamic programming sequential approach. To address the uncertainty of departure/arrival times in recovery plans, the authors divided all flights into two parts: the flights with deterministic departure/arrival times and those with several stochastic departure/arrival times. \citeauthor{maharjan2011optimization} (\citeyear{maharjan2011optimization}) developed a binary quadratic integer programming model to recover gate assignments and heuristically solved the model by a gate zoning strategy. \citeauthor{zhang2017optimization} (\citeyear{zhang2017optimization}) built two network flow models for the gate reassignment problem without/with connecting passengers. To overcome the difficulty in solving large-scale problems, the authors used a gate clique technique and a variable rolling horizon strategy. \citeauthor{poyraz2024} (\citeyear{poyraz2024}) built a bi-objective gate reassignment model to balance the efficiency (using the gates as much as possible) and stability (sticking as closely as possible to the initial plan). To improve the solution efficiency, they decomposed their model into small-scale models by clustering the aircraft. In this paper, we extend the schedule and aircraft recovery problem by further integrating the gate reassignment problem. Then we develop a Benders and column generation method to solve the above integrated recovery problem.\par

\subsection{Integrated recovery}

\citeauthor{su2021airline} (\citeyear{su2021airline}) provided a state-of-the-art review, focusing on the characteristics of mathematical models and solution methods in disruption management. They suggested that airline disruption management should focus on integrating recovery stages and reducing the complexity of the integrated models. Motivated by this work, we consider the integrated airline recovery problem in this paper and apply some effective techniques to alleviate the solution difficulty by exploring the specific structure.\par

\textit{Airline integrated recovery problem} seeks to preserve the interdependence among multiple recovery phases in airline disruption recovery. Various combinations of recovery phases have been studied and different methods have been developed to solve the integrated airline recovery problem (e.g., \citeauthor{bratu2006flight}, \citeyear{bratu2006flight}; \citeauthor{petersen2012optimization}, \citeyear{petersen2012optimization}; \citeauthor{sinclair2014}, \citeyear{sinclair2014}, \citeyear{sinclair2016}; \citeauthor{maher2015novel}, \citeyear{maher2015novel}, \citeyear{maher2016solving}; \citeauthor{arikan2016integrated}, \citeyear{arikan2016integrated}, \citeyear{arikan2017flight};  \citeauthor{hu2016integrated}, \citeyear{hu2016integrated}; \citeauthor{zhang2016math}, \citeyear{zhang2016math}; \citeauthor{zhang2017optimization}, \citeyear{zhang2017optimization}; \citeauthor{pternea2019aircraft}, \citeyear{pternea2019aircraft}; \citeauthor{hu2021integrated}, \citeyear{hu2021integrated}; \citeauthor{cadarso2022}, \citeyear{cadarso2022}; \citeauthor{ding2023towards}, \citeyear{ding2023towards};  \citeauthor{xu2023distributionally}, \citeyear{xu2023distributionally}; \citeauthor{zhong2024}, \citeyear{zhong2024}; \citeauthor{JIANG2025104243}, \citeyear{JIANG2025104243};
\citeauthor{wang2025}, \citeyear{wang2025}; etc.).\par

\citeauthor{petersen2012optimization} (\citeyear{petersen2012optimization}) formulated a fully integrated airline recovery problem that included flight, aircraft, crew, and passenger recovery phases by combining several path-based models. The authors limited the recovery scope by identifying disruptable flights and employing event-driven delays. A Benders decomposition framework was proposed to link these models, where new routings and pairings were generated by the CG method. The authors proposed effective certificates to judge the infeasibility and sub-optimality of the Benders subproblem with linear equality constraints during the process of CG. Motivated by it, we extend their infeasibility certificates to further consider the capacity inequality constraints in Benders subproblem.\par

\citeauthor{maher2015novel} (\citeyear{maher2015novel}) considered the fully integrated airline recovery problem and developed a path-based model to solve it. The author used the column-and-row generation method that was proposed by \citeauthor{maher2016solving} (\citeyear{maher2016solving}) to solve the integrated model and provided a thorough analysis. Unlike \citeauthor{maher2015novel} (\citeyear{maher2015novel}, \citeyear{maher2016solving}), we use the Benders decomposition framework, CG method, and reformulation techniques to reduce the complexity of our integrated recovery problem by dividing the integrated model into several small-scale problems.\par

\citeauthor{arikan2017flight} (\citeyear{arikan2017flight}) considered a fully integrated recovery model by taking crew recovery into account in \citeauthor{arikan2016integrated} (\citeyear{arikan2016integrated}). In addition to regular recovery options, they employed cruise speed control. The authors developed a network flow model and then reformulated it as a conic quadratic mixed integer programming that can be tackled by commercial solvers. They applied approximate flight delay and eliminated unlikely used connection arcs to tackle the large airline network. Benefiting from the favorable structure in our integrated problem, we propose some techniques to reduce the scale of our Benders subproblem while the associated recovery performance is not sacrificed (see Section \ref{Section Benders subproblem (BSP)} and Section \ref{Separation} for details).\par

Based on the fully integrated recovery model in \citeauthor{arikan2017flight} (\citeyear{arikan2017flight}),  \citeauthor{ding2023towards} (\citeyear{ding2023towards}) developed a mathematical model considering more details about crew and passenger recovery. The authors proposed a solution framework by incorporating deep reinforcement learning into the variable neighborhood search algorithm, which can quickly obtain an integrated recovery plan. Since this strategy cannot provide optimality gaps for their solutions, the solution quality of the associated plan might not be stable. In our work, we focus on developing a method to find high-quality solutions with optimality gaps within a reasonable runtime.

In summary, airline integrated recovery has received significant attention in the literature. From the perspective of an integrated recovery paradigm, further integrating airport gate reassignment into the schedule and aircraft recovery problem has the potential to minimize disruptions to airport resource allocation and reduce recovery costs. To the best of our knowledge, this is the first work in airline disruption management to explore the benefits of such an integration.\par

\section{Airline integrated recovery with gate reassignment}\label{Section2}
In this section, we first provide the definitions of several key terms used in this paper. Then, we develop a mathematical model to formulate the airline integrated recovery problem.\par

\subsection{Definitions}\label{Section Definitions}
\begin{figure}[htbp]
   \hspace{-0.14cm}
   \begin{minipage}{0.46\linewidth}
    \centering
    {\includegraphics[width=0.85\textwidth]{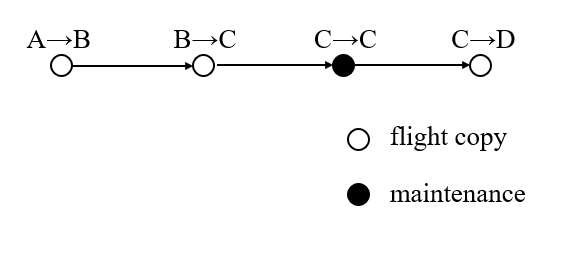}}
    \caption{Illustration of an aircraft route\label{aircraftRoute}}
   \end{minipage}
   \qquad
   \hspace{-1.3cm}
   \begin{minipage}{0.46\linewidth}
    \centering
    {\includegraphics[width=1.18\textwidth]{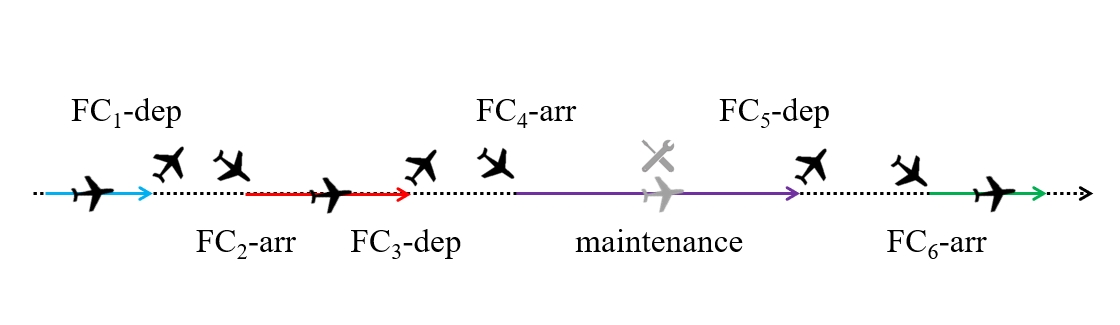}}
    \vspace*{-0.25cm}
    \caption{Illustration of a gate pattern\label{gatePattern}}
   \end{minipage}
 \end{figure}

\begin{itemize}
\setlength{\itemsep}{0pt}
\setlength{\parsep}{0pt}
\setlength{\parskip}{0pt}

 \item[1)] \textit{Recovery time window}: The Recovery time window is a time period during which all recovery actions are taken. Airline operation is expected to return to its normal state before the end of the recovery time window. The length of the recovery time window mainly depends on the severity of the disruption and the resource limits of stakeholders.

 \item[2)] \textit{Recovery options}: There are various recovery options available in disruption management. This paper uses three common recovery options as candidates: flight delays, flight cancellations, and swapping flights. To simplify the integrated recovery problem, we assume that aircraft swapping and gate swapping (i.e., changing the gate assignments for flights) are restricted to the same types of equipment, but the SAGRM formulation is also applicable to the case where flights are swapped to different types of equipment.

  \item[3)] \textit{Slot}: A slot is a time period that has the following attributes: airport, start time, time interval, and slot capacity. The aircraft movement capacity of a slot at a particular airport might be decreased once a disruption occurs.

  \item[4)]\textit{Flight copy}: We use the uniform flight delay copy technique to model flight delays. That is, given a maximum allowable delay, \textit{copies} of an original flight at several uniform time intervals are made to represent candidate delayed flights. The copy of the original flight is called a \textit{flight copy} in this paper. In addition, the flying time of a flight copy is assumed to be the same as that of the original flight, but it can be easily considered in the SAGRM formulation by introducing the flight copies with different flying times.

  \item[5)]\textit{Aircraft route}: A route of an aircraft is a sequence of flight copies performed by the aircraft, which is illustrated in Figure \ref{aircraftRoute}. Two letters on each circle dot in Figure \ref{aircraftRoute} denote the associated departure and arrival airports, respectively. The cost of an aircraft route in the disruption recovery context is related to the total delay time and the number of swapped flights in the route. \\
      \hspace*{0.4cm}Some key constraints in forming an aircraft route are presented as follows. (a) Type constraint: The route should only contain the flight copies whose types are consistent with the type of the associated aircraft. (b) Connection constraint: In addition to the space and time matches between each two consecutive flight copies in the route, the start and end airports of a route should be the airports where the aircraft is located before and after the recovery time window, respectively. (c) Maintenance constraint: The planned maintenance requirement can be regarded as a \textit{special flight} that a route must contain for a specific aircraft.

  \item[6)]\textit{Gate pattern}: A gate pattern for a particular gate type (e.g., wide-body, narrow-body, close-in or remote) at a particular airport is a sequence of flight copy activities (representing departure, arrival, and maintenance) which is serviced by a particular gate; see Figure \ref{gatePattern} for an illustration. In this paper, the cost of one pattern depends on the location of the airport and the associated gate type. In the existing literature, the cost of gate reassignment is also related to the time and space inconsistency (\citeauthor{yan2011airport}, \citeyear{yan2011airport}). Time inconsistency has been considered in the total delay time of aircraft routes in our work. Since the cost caused by space inconsistency is relatively small during disruption management in practice, this cost is omitted in this paper.\\ \label{defination gate patterns}
      \hspace*{0.4cm}Two key constraints in forming a gate pattern are presented as follows. (a) Match constraint: A pattern only contains the flight copy activities whose types are consistent with the type of the associated gate. (b) Buffer time constraint: The time gap between a departure flight copy and the next upcoming arrival flight copy at a particular gate in a gate pattern should be more than a specified buffer time.
 \end{itemize}

\subsection{Mathematical model}\label{Section Formulation 1}
We now develop a path-based schedule, aircraft, and gate recovery model (SAGRM) to formulate the integrated recovery problem for flight schedule, aircraft, and gate. First, we give the notations used in this paper as follows. Then, the SAGRM formulation is presented.\par

\begin{center}
\renewcommand{\baselinestretch}{1.1}
\begin{small}
	\begin{longtable}{p{0.15\textwidth}p{0.75\textwidth}}
    \caption{Notations used in this paper}\label{Notations}
		\endfirsthead
			
        \midrule
		\endhead
			
        \midrule
        \multicolumn{2}{r}{\footnotesize ({\it continued on next page})}
		\endfoot
			
		\hline
	\endlastfoot

    \toprule
    \multicolumn{2}{l}{\bf{Sets}}\\
    \midrule
    $S$ &the set of slots\\
    $F$ &the set of original flights\\
    $FC$ &the set of flight copies\\
    $FC_j$ &the set of flight copies for original flight~$j\in F$\\
    $E$ &the set of feasible connections between two flight copies\\
    $R$ &the set of aircraft\\
    $P^r$ &the set of feasible routes for aircraft $r\in R$\\
    $\bar{P}^r$ &the subset of feasible routes for aircraft $r\in R$, i.e., $\bar{P}^r \subseteqq P^r$ \\
    $A$ &the set of airports\\
    $T_{a}$ &the set of gate types at airport $a\in A$\\
    $P^{a,t}$ &the set of patterns in airport $a\in A$ with gate type $t\in T_{a}$\\
    $\bar{P}^{a,t}$ &the subset of patterns in airport $a\in A$ with gate type $t\in T_{a}$, i.e., $\bar{P}^{a,t} \subseteq P^{a,t}$ \\
    $FC^{a,t}_{dep}$ &the set of flight copies that depart at airport $a\in A$ with gate type $t\in T_{a}$\\
    $FC^{a,t}_{arr}$ &the set of flight copies that arrive at airport $a\in A$ with gate type $t\in T_{a}$\\
    $E^{a,t}$ &the set of feasible connections at airport $a\in A$ with gate type $t\in T_{a}$\\
    \toprule
    \multicolumn{2}{l}{\bf{Parameters}}\\
    \midrule
    $u_{sdep}/u_{sarr}$ &the departure/arrival capacity of slot $s\in S$\\
    $d_j$ &the cost of canceling original flight $j\in F$\\
    $c_p^r$ &the cost of aircraft route $p\in P^r$\\
    $a_{{p,{j_v}}}^{r}$ &1, if flight copy $j_v\in FC_j$ is in route $p\in P^r$;~0, otherwise\\
    $a_{{p,{j_v}}}^{r,sdep}/a_{{p,{j_v}}}^{r,sarr}$ &1, if flight copy $j_v\in FC_j$ is in route $p\in P^r$ and departs/arrives in slot $s \in S$; \\
    &0, otherwise\\
    $b_{p,(i_u,j_v)}^{r}$ &1, if the connection between flight copy $(i_u,j_v) \in E$, ~is in route $p\in P^r$;~0, otherwise\\
    $c_{a,t}$ &the cost of a gate pattern at airport $a\in A$ with gate type $t\in T_{a}$\\
    $n_{a,t}$ &the total number of gates at airport $a\in A$ with gate type $t\in T_{a}$\\
    $a_{{p,{j_{v}dep}}}^{a,t}/a_{{p,{j_{v}arr}}}^{a,t}$ &1, if departure/arrival activity of flight copy $j_v\in FC_j$ is in pattern $p\in P^{a,t}$;\\
    &0, otherwise\\
    $b_{p,(i_u,j_v)}^{a,t}$ &1, if connection $(i_u,j_v) \in E$~is in pattern $p\in P^{a,t}$;~0, otherwise\\
    \toprule
    \multicolumn{2}{l}{\bf{Decision variables}}\\
    \midrule
    $z_j$ &1, if original flight~$j\in F$ is canceled;~0, otherwise\\
    $y_p^r$ &1, if aircraft $r\in R$ is assigned to route $p\in P^r$;~0, otherwise\\
    $w_p^{a,t}$ &1, if gate pattern $p\in P^{a,t}$ is chosen;~0, otherwise\\
    \toprule
    \multicolumn{2}{l}{\bf{Abbreviations}}\\
    \midrule
    SAGRM & Schedule, Aircraft, and Gate Recovery Model (\ref{SAGobj})$-$(\ref{SAGcon11}) \\
    SARM & Schedule and Aircraft Recovery Model (\ref{SARMobj})$-$(\ref{SARMcon9}) \\
    GRM & Gate Recovery Model (\ref{GRMobj})$-$(\ref{GRMcon5}) \\
    LR-X & the linear relaxation of (mixed) integer problem X\\
    CG & column generation\\
    BMP & Benders master problem\\
    RxBMP & the relaxed BMP, i.e., the BMP with the partial Benders cuts\\
    RtRxBMP & the restricted RxBMP, i.e., the RxBMP with partial routes $\bar{P}^{r}$\\
    BSP & Benders subproblem\\
    RtBSP & the restricted BSP, i.e., the BSP with partial gate patterns $\bar{P}^{a,t}$\\
    P1LR-RtBSP & the phase one optimization problem of LR-RtBSP\\
    \bottomrule
	\end{longtable}
\end{small}
\end{center}

As mentioned in the motivations for integrating three recovery phases in Section \ref{Introduction}, the sequential way of recovering these phases may lead to an undesirable result, especially in an airline network. In order to preserve the interdependence among these recovery phases, we propose the SAGRM formulation that integrates the recovery of flight schedule, aircraft route, and gate assignment. Schedule recovery and aircraft rerouting are crucial in disruption management since their decisions could directly affect other recovery phases (\citeauthor{XIONG201364}, \citeyear{XIONG201364}; \citeauthor{WU201962}, \citeyear{WU201962}; \citeauthor{BRUECKNER2021102333}, \citeyear{BRUECKNER2021102333}). These two phases are often considered simultaneously to obtain a high-quality recovery plan. To reduce the complexity of the integrated schedule and aircraft recovery model, some work uses a path-based modelling approach to construct candidate routes (\citeauthor{petersen2012optimization}, \citeyear{petersen2012optimization}; \citeauthor{maher2016solving}, \citeyear{maher2016solving}; \citeauthor{liang2018column}, \citeyear{liang2018column}). For the gate reassignment problem, the previous work usually focuses on a single airport and uses an arc-based model to describe all practical restrictions (\citeauthor{yan2009airport}, \citeyear{yan2009airport}; \citeauthor{yan2011airport}, \citeyear{yan2011airport}; \citeauthor{zhang2017optimization}, \citeyear{zhang2017optimization}). On the one hand, due to the network effect in an airline network, the recovery decision of gate reassignment in an airport would also affect other airports¡¯ decisions. To capture this effect, we formulate the integrated recovery at multiple airports. On the other hand, the arc-based modelling approach requires a lot of variables and constraints to formulate the problem. The SAGRM uses the path-based modelling approach to describe aircraft routes and gate patterns. Benefiting from it, we can easily establish the relationship among three recovery phases in an airline network.\par

The specific formulation of SAGRM consists of three parts corresponding to the integrated schedule and aircraft recovery, the gate reassignment, and the linking constraints. The path-based formulation for the first part is based on the previous literature (\citeauthor{petersen2012optimization}, \citeyear{petersen2012optimization}; \citeauthor{liang2018column}, \citeyear{liang2018column}). Motivated by the gate assignment model in \citeauthor{diepen2012gate} (\citeyear{diepen2012gate}), we develop a path-based formulation for the gate reassignment problem in SAGRM. According to the structure of the integrated recovery problem, the linking constraints are introduced to link the above two parts of formulations. Using the notations in Table \ref{Notations}, we now provide the mathematical formulation of the SAGRM as follows.\par

\vspace{0.5em}
\textbf{SAGRM Formulation:}\par
\vspace{-1.5em}
\begin{eqnarray}
&\min &\sum_{j\in {F}}d_j z_j +\sum_{r\in R}\sum_{p\in {P^r}}c_p^r y_p^r
        +\sum_{a\in A}\sum_{t\in T_{a}}\sum_{p\in {P^{a,t}}}c_{a,t} w_p^{a,t}\label{SAGobj}\\
&\mathrm{s.t.} &\sum_{r\in R}\sum_{p\in {P^r}}\sum_{j_v\in {FC_j}}a_{{p,{j_v}}}^{r} y_p^r+z_{j}=1,\ \ \forall j\in {F},\label{SAGcon1}\\
&&
\sum_{r\in {R}}\sum_{p\in {P^r}}\sum_{j\in {F}}\sum_{j_v\in {FC_j}} a_{{p,{j_v}}}^{r,sdep} y_p^r \leq u_{sdep},\ \ \forall s\in S,\label{SAGcon2}\\
&&
\sum_{r\in {R}}\sum_{p\in {P^r}}\sum_{j\in {F}}\sum_{j_v\in {FC_j}} a_{{p,{j_v}}}^{r,sarr} y_p^r \leq u_{sarr},\ \ \forall s\in S,\label{SAGcon3}\\
&&
\sum_{p\in {P^r}} y_p^r \leq 1,\ \ \forall r\in R,\label{SAGcon4}\\
&&
\sum_{a\in A}\sum_{t\in T_{a}}\sum_{p\in {P^{a,t}}}a_{{p,{j_{v}dep}}}^{a,t} w_p^{a,t}
     = \sum_{r\in R}\sum_{p\in {P^r}}a_{{{p},j_v}}^{r} y_p^r,\ \ \forall j\in {F},{j_v\in {FC_j}}, \label{SAGcon5}\\
&&
\sum_{a\in A}\sum_{t\in T_{a}}\sum_{p\in {P^{a,t}}}a_{{p,{j_{v}arr}}}^{a,t} w_p^{a,t}
     = \sum_{r\in R}\sum_{p\in {P^r}}a_{{{p},j_v}}^{r} y_p^r,\ \ \forall j\in {F},{j_v\in {FC_j}}, \label{SAGcon6}\\
&&
\sum_{a\in A}\sum_{t\in T_{a}}\sum_{p\in {P^{a,t}}}b_{{p,{(i_u,j_v)}}}^{a,t} w_p^{a,t}
     \geq \sum_{r\in R}\sum_{p\in {P^r}}b_{p,(i_u,j_v)}^{r} y_p^r ,\ \ \forall (i_u,j_v)\in {E}, \label{SAGcon7}\\
&&
\sum_{p\in {P^{a,t}}} w_p^{a,t} \leq n_{a,t},\ \ \forall a\in A, t\in T_{a},\label{SAGcon8}\\
&&
z_{j}\in \{0,1\},\ \ \forall j\in {F}, \label{SAGcon9}\\
&&
y_p^r\in \{0,1\},\ \ \forall r\in R, p\in P^r,\label{SAGcon10}\\
&&
w_p^{a,t}\in \{0,1\},\ \ \forall a\in A, t\in T_{a}, p\in P^{a,t}.\label{SAGcon11}
\end{eqnarray}

The objective (\ref{SAGobj}) is to minimize the total recovery cost from flight cancellation, aircraft rerouting, and gate reassignment. The calculation of the cost of routes and gate patterns has been discussed earlier in Section \ref{Section Definitions}. Constraints (\ref{SAGcon1}) ensure each original flight $j$ is either canceled or contained in one chosen route. Departure and arrival capacity restrictions for slots are captured in (\ref{SAGcon2}) and (\ref{SAGcon3}), respectively. Constraints (\ref{SAGcon4}) require that each aircraft $r$ is assigned to exactly one route or is idle during the disruption period. Constraints (\ref{SAGcon5}) and (\ref{SAGcon6}) ensure the recovery option for each flight must be consistent among three recovery phases (flight schedule, aircraft, and gate). For each original flight $j$, the chosen flight copy activities $j_{v}dep$ and $j_{v}arr$ in the gate recovery should be consistent with the chosen flight copy $j_{v}$ in the aircraft recovery. For any two consecutive flight copies $(i_u,j_v)$, Constraints (\ref{SAGcon7}) guarantee that the arrival gate of $i_u$ and the departure gate of $j_v$ are the same if $(i_u,j_v)$ are contained in a chosen route by the same aircraft. Constraints (\ref{SAGcon8}) make sure that the number of chosen gate patterns does not exceed the number of available gates $n_{a,t}$. In addition, the key constraints in forming aircraft routes and gate patterns, as mentioned in Section \ref{Section Definitions} 5)-6), are considered when constructing sets $P^r$ and $P^{a,t}$ respectively, which are provided in \ref{appendix CP models}.\par

It is noted that the gate recovery solution derived from SAGRM corresponds to different gate types at airports {\it{rather than to specific gates}}. The main benefit of this modelling approach is that the scale of the model is greatly reduced. To obtain the actual reassignment for each gate, the only thing to do is to allocate the selected gate patterns (the gate patterns whose $w_p^{a,t}$ are equal to 1 in SAGRM) to the gates according to their located airports $a$ and types $t$. In the following sections, we will utilize the structure of the model to design an effective algorithm for SAGRM and propose acceleration techniques for the algorithm by exploiting the characteristics of the integrated recovery problem.\par

\section{Solution methodology}\label{Section3}
SAGRM is a large-scale integer programming problem with complex constraints, so it is not expected that an optimal solution can be obtained efficiently through commercial solvers within a reasonable runtime. Due to the inherent multi-stage structure of SAGRM, it allows us to utilize a Benders decomposition framework to solve SAGRM effectively by dividing the integrated model into two smaller sub-models. Accordingly, the Benders master problem (BMP) of SAGRM contains two recovery phases: flight rescheduling and aircraft rerouting. Given a schedule and aircraft recovery solution, the Benders subproblem (BSP) of SAGRM considers gate reassignment and generates a Benders cut so as to adjust the current schedule and aircraft recovery solution. To avoid enumerating all aircraft routes and gate patterns, a column generation technique is used to iteratively generate the routes in $P^r$ and the patterns in $P^{a,t}$. Thus, we propose a solution method that combines Benders decomposition and column generation (BCG) to solve SAGRM.\par

\subsection{Benders subproblem (BSP)}\label{Section Benders subproblem (BSP)}
Let $(\textbf{y}, \textbf{z})$ be the set of the solutions that satisfy Constraints (\ref{SAGcon1})$-$(\ref{SAGcon4}) and (\ref{SAGcon10}). For any given solution $\bar{y}_p^r \in \textbf{y}$, the associated BSP, which is a gate recovery model (GRM) corresponding to $\bar{y}_p^r$, can be obtained from SAGRM as follows.\par

\vspace{0.2em}
\textbf{GRM Formulation as BSP:}
\begin{eqnarray}
&\min &\sum_{a\in A}\sum_{t\in T_{a}}\sum_{p\in {P^{a,t}}}c_{a,t} w_p^{a,t} \label{GRMobj}\\
&\mathrm{s.t.} &\sum_{a\in A}\sum_{t\in T_{a}}\sum_{p\in {P^{a,t}}}a_{{p,{j_{v}dep}}}^{a,t} w_p^{a,t}
     = \sum_{r\in R}\sum_{p\in {P^r}}a_{{{p},j_v}}^{r} \bar{y}_p^r,\ \ \forall j\in {F},{j_v\in {FC_j}},\label{GRMcon1}\\
&&
     \sum_{a\in A}\sum_{t\in T_{a}}\sum_{p\in {P^{a,t}}}a_{{p,{j_{v}arr}}}^{a,t} w_p^{a,t}
     = \sum_{r\in R}\sum_{p\in {P^r}}a_{{{p},j_v}}^{r} \bar{y}_p^r,\ \ \forall j\in {F},{j_v\in {FC_j}},\label{GRMcon2}\\
&&
     \sum_{a\in A}\sum_{t\in T_{a}}\sum_{p\in {P^{a,t}}}b_{{p,{(i_u,j_v)}}}^{a,t} w_p^{a,t}
     \geq \sum_{r\in R}\sum_{p\in {P^r}}b_{p,(i_u,j_v)}^{r} \bar{y}_p^r ,\ \ \forall (i_u,j_v)\in {E}, \label{GRMcon3}\\
&&
     \sum_{p\in {P^{a,t}}} w_p^{a,t} \leq n_{a,t},\ \ \forall a\in A, t\in T_{a}, \label{GRMcon4}\\
&&
     w_p^{a,t}\in \{0,1\},\ \ \forall a\in A, t\in T_{a}, p\in P^{a,t}.\label{GRMcon5}
\end{eqnarray}

In the standard Benders decomposition framework, the Benders subproblem is a linear programming (LP) problem. However, such a decomposition strategy can also be taken to solve pure integer programming problems.  In such a case, the associated Benders subproblem is an integer linear programming problem and then the Benders cuts are obtained by its linear relaxation (\citeauthor{wolsey2020integer}, \citeyear{wolsey2020integer}). The above manner has been widely used to deal with the problems in the airline industry (e.g., \citeauthor{cordeau2001benders}, \citeyear{cordeau2001benders}; \citeauthor{papadakos2009integrated}, \citeyear{papadakos2009integrated}; \citeauthor{petersen2012optimization}, \citeyear{petersen2012optimization}; etc.). To get the Benders cuts, we consider the linear relaxation of the BSP (LR-BSP). Let $\pi_{j_{v}dep}, \pi_{j_{v}arr}, \pi_{(i_u,j_v)}, \pi_{a,t}$ be the dual values related to Constraints (\ref{GRMcon1})$-$(\ref{GRMcon4}) respectively, and $\pi$ be the vector $(\pi_{j_{v}dep}, \pi_{j_{v}arr}, \pi_{(i_u,j_v)})$. The representation of Benders feasibility and optimality cuts are:

\vspace*{-0.5cm}
\begin{align}
&\sum_{r\in R}\sum_{p\in {P^r}}H_{p}^{r}(\pi)y_p^r - \sum_{a\in A}\sum_{t\in T_{a}} \pi_{a,t}n_{a,t}\leq 0,
\ \ \forall (\pi,\pi_{a,t}) \in \Pi_{Fea}, \label{classFeaCut}\\
&\sum_{r\in R}\sum_{p\in {P^r}}H_{p}^{r}(\pi)y_p^r - \sum_{a\in A}\sum_{t\in T_{a}} \pi_{a,t}n_{a,t}\leq q,
\ \ \forall (\pi,\pi_{a,t})\in \Pi_{Opt},\label{classOptCut}
\end{align}

\noindent where $H_{p}^{r}(\pi)= \sum\limits_{j\in {F}}\sum\limits_{j_v\in {FC_j}}(\pi_{j_{v}dep}+\pi_{j_{v}arr})a_{{{p},j_v}}^{r} + \sum\limits_{(i_u,j_v)\in {E}}\pi_{(i_u,j_v)}b_{p,(i_u,j_v)}^{r}$. The sets $\Pi_{Fea}$ and $\Pi_{Opt}$ denote the set of extreme points and extreme directions of the dual feasible region of the LR-BSP, respectively, and $q$ is a decision variable linking BMP and BSP.\par

As mentioned before, CG is applied to deal with the large size of $P^{a,t}$, which iterates between the CG-restricted master problem and CG-subproblem. During the iteration of CG, the gate patterns set $P^{a,t}$ is restricted to $\bar{P}^{a,t} \subseteq P^{a,t}$, i.e., the BSP becomes the restricted BSP (RtBSP). The associated CG-restricted master problem is the linear relaxation of the RtBSP. Given its dual solutions, the CG-subproblem is a resource-constrained shortest path problem that finds the patterns with the most negative reduced cost on a flight connection network. This network is similar to the aircraft connection network in \citeauthor{papadakos2009integrated} (\citeyear{papadakos2009integrated}), but the nodes in our network represent flight copy activities. Given airport $a$ and gate type $t$, the reduced cost $\bar{c}_p^{a,t}$ of gate pattern $p \in P^{a,t}$ is

\vspace*{-0.5cm}
\begin{align}
 \bar{c}_p^{a,t} :=
  c_{a,t}- \sum_{j\in {F}}\sum_{j_v\in {FC_j}}(\pi_{j_{v}dep}a_{{p,{j_{v}dep}}}^{a,t}+\pi_{j_{v}arr}a_{{p,{j_{v}arr}}}^{a,t})
 - \sum_{(i_u,j_v)\in {E}}\pi_{(i_u,j_v)}b_{p,(i_u,j_v)}^{a,t} + \pi_{a,t}. \label{patternReducedCost}
\end{align}

\noindent CG-subproblem typically pursues a gate pattern with the most negative reduced cost. Besides, multiple patterns with negative reduced costs could be added to $\bar{P}^{a,t}$ to accelerate the CG.

It is noted that the information from solution $\bar{y}_p^r$ can be exploited to improve the efficiency of our BCG method. The following three well-known techniques are used in this paper.
\begin{itemize}
\setlength{\itemsep}{3pt}
\setlength{\parsep}{0pt}
\setlength{\parskip}{0pt}

 \item[(1)] {\it Removal of redundant constraints}: The number of Constraints (\ref{GRMcon3}), equal to the cardinality of the feasible connections set $E$, is much larger than the number of other constraints in the BSP. Fortunately, lots of these constraints can be reduced for a given RxBMP solution $\bar{y}_{p}^{r}$. Recall that the constraints in (\ref{GRMcon3}) are redundant when their right-hand terms $\sum\limits_{r\in R}\sum\limits_{p\in {P^r}}b_{p,(i_u,j_v)}^{r} \bar{y}_{p}^{r}$ are equal to $0$. Hence, we can remove these redundant constraints by using the information from $\bar{y}_{p}^{r}$.

\item[(2)] {\it Customized initialization of gate patterns}: Based on the information in $\bar{y}_{p}^{r}$, we can accelerate the CG process by generating some gate patterns that are more likely to be selected before using CG to solve the BSP, which is implemented as follows. Suppose flight copy connection $(i_u, j_v)$ is in the chosen routes related to $\bar{y}_{p}^{r}$. If $(i_u, j_v)$ has not been covered in the previous iterations, we will generate several gate patterns containing the arrival of $i_u$ and the departure of $j_v$. This task can be carried out by constraint programming solvers, e.g., CPLEX.\par

\item[(3)] {\it Application of valid inequalities}: The Benders cuts (\ref{classFeaCut})$-$(\ref{classOptCut}) are obtained by the dual of LR-BSP so that these cuts may be weak for BMP due to the integrality gap between BSP and LR-BSP. It is noted that the linking variables $\bar{y}_p^r$ are all binary, then some specific valid inequalities can be applied to tackle the gap. In this paper, we apply the following three valid inequalities: no-good cut (\citeauthor{wolsey2020integer}, \citeyear{wolsey2020integer}), Laporte \& Louveaux cut (\citeauthor{laporte1993integer}, \citeyear{laporte1993integer}), and global cut (\citeauthor{fakhri2017benders}, \citeyear{fakhri2017benders}).\par
\end{itemize}\par

In Section \hyperref[Section4]{5}, we will further exploit the special structure of the BSP and the known information from solution $\bar{y}_p^r$ to develop acceleration techniques for our BCG method.\par

\subsection{Benders master problem (BMP)}
The BMP, a schedule and aircraft recovery model (SARM), is formulated by replacing Constraints (\ref{SAGcon5})$-$(\ref{SAGcon8}) and (\ref{SAGcon11}) in the SAGRM with the Benders cuts (\ref{classFeaCut}) and (\ref{classOptCut}). In addition, the gate recovery cost $\sum\limits_{a\in A}\sum\limits_{t\in T_{a}}\sum\limits_{p\in {P^{a,t}}}c_{a,t} w_p^{a,t}$ in objective (\ref{SAGobj}) is replaced with the decision variable $q$ in the Benders optimality cuts (\ref{classOptCut}). Thus, the SARM formulation as BMP is as follows.

\vspace{0.5em}
\textbf{SARM Formulation as BMP:}
\begin{eqnarray}
&\min&\sum_{j\in {F}}d_j z_j +\sum_{r\in R}\sum_{p\in {P^r}}c_p^r y_p^r + q \label{SARMobj}\\
&\mathrm{s.t.}&\sum_{r\in R}\sum_{p\in {P^r}}\sum_{j_v\in {FC_j}}a_{{p,{j_v}}}^{r} y_p^r+z_{j}=1,\ \ \forall j\in {F}, \label{SARMcon1}\\
&&
      \sum_{r\in {R}}\sum_{p\in {P^r}}\sum_{j\in {F}}\sum_{j_v\in {FC_j}} a_{{p,{j_v}}}^{r,sdep} y_p^r \leq u_{sdep},\ \ \forall s\in S, \label{SARMcon2}\\
&&
      \sum_{r\in {R}}\sum_{p\in {P^r}}\sum_{j\in {F}}\sum_{j_v\in {FC_j}} a_{{p,{j_v}}}^{r,sarr} y_p^r \leq u_{sarr},\ \ \forall s\in S, \label{SARMcon3}\\
&&
     \sum_{p\in {P^r}} y_p^r \leq 1,\ \ \forall r\in R,
     \label{SARMcon4}\\
&&
     \sum_{r\in R}\sum_{p\in {P^r}}H_{p}^{r}(\pi)y_p^r - \sum_{a\in A}\sum_{t\in T_{a}} \pi_{a,t}n_{a,t}\leq 0,
     \ \ \forall (\pi,\pi_{a,t}) \in \Pi_{Fea}, \label{SARMcon5}\\
&&
     \sum_{r\in R}\sum_{p\in {P^r}}H_{p}^{r}(\pi)y_p^r - \sum_{a\in A}\sum_{t\in T_{a}} \pi_{a,t}n_{a,t}\leq q,
     \ \ \forall (\pi,\pi_{a,t})\in \Pi_{Opt}, \label{SARMcon6}\\
&&
     z_{j}\in \{0,1\},\ \ \forall j\in {F}, \label{SARMcon7}\\
&&
     y_p^r\in \{0,1\},\ \ \forall r\in R, p\in P^r \label{SARMcon8}\\
&&
     q\geq0 \label{SARMcon9}.
\end{eqnarray}

The Benders cuts (\ref{SARMcon5}) and (\ref{SARMcon6}) are valid because of the gap between the BSP and its linear relaxation. Due to the large size of $\Pi_{Fea}$ and $\Pi_{Opt}$, the elements in two sets are gradually obtained during the iteration process of the BCG method. Hence, Constraints (\ref{SARMcon5}) and (\ref{SARMcon6}) are relaxed to the constraints related to $\bar{\Pi}_{Fea} \subseteq \Pi_{Fea}$ and $\bar{\Pi}_{Opt} \subseteq \Pi_{Opt}$ respectively in the iteration, and the BMP turns into a relaxed Benders master problem (RxBMP).\par

Then, due to the large size of aircraft routes in $P^r$, the RxBMP is solved by the CG technique. During the iteration process of CG, the RxBMP turns into the restricted RxBMP (RtRxBMP), i.e., the aircraft routes set $P^{r}$ is restricted to $\bar{P}^{r} \subseteq P^{r}$. The CG-restricted master problem for the RxBMP is the linear relaxation of the RtRxBMP. The associated CG-subproblem turns into a resource-constrained shortest path problem on a flight copy connection network where the nodes represent flight copies. Given aircraft $r$, the reduced cost $\bar{c}_p^{r}$ corresponding to the aircraft route $p \in P^r$ is:
\vspace{-0.2cm}

\begin{align}
 \nonumber
  \bar{c}_p^{r} :=\ &c_p^{r}- \sum_{j\in {F}}\sum_{j_v\in {FC_j}}\lambda_{j} a_{{p,{j_v}}}^{r}
 + \sum_{s\in S}\sum_{j\in {F}}\sum_{j_v\in {FC_j}}(\lambda_{sdep}a_{{p,{j_v}}}^{r,sdep} + \lambda_{sarr}a_{{p,{j_v}}}^{r,sarr}) + \lambda_{r} \\
 &
  ~ + \sum_{\pi\in {\bar{\Pi}_{fea}}} \lambda_{\pi}^{fea}H_{p}^{r}(\pi) +  \sum_{\pi\in {\bar{\Pi}_{opt}}} \lambda_{\pi}^{opt}H_{p}^{r}(\pi),
  \label{routeReducedCost}
\end{align}

\noindent
where $\lambda_{j}, \lambda_{sdep}, \lambda_{sarr}, \lambda_{r}, \lambda_{\pi}^{fea}$ and $\lambda_{\pi}^{opt}$ are the dual values of Constraints (\ref{SARMcon1})$-$(\ref{SARMcon6}), and $\bar{\Pi}_{fea}$/$\bar{\Pi}_{opt}$ denotes the set of vectors $(\pi_{j_{v}dep}, \pi_{j_{v}arr}, \pi_{(i_u,j_v)})$. Similar to the way of applying the CG technique to the BSP, we can seek an appropriate number of additional aircraft routes with negative reduced costs to decrease the computational time.\par

\subsection{BCG method}
\begin{figure}[htbp]
\centering
\includegraphics[width=0.63\textwidth]{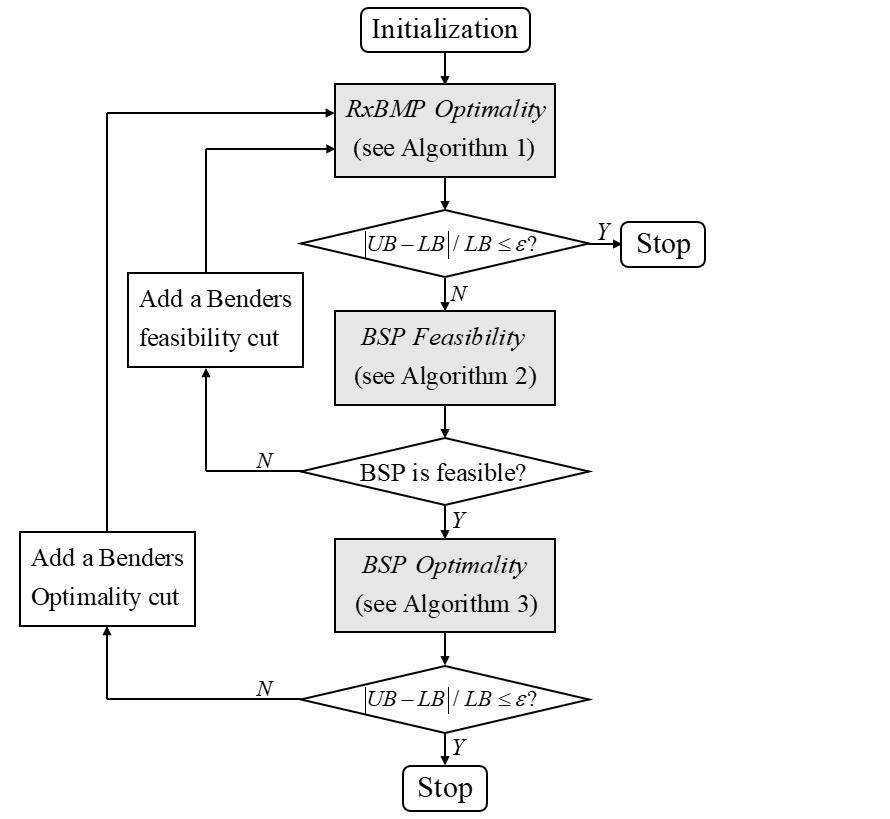}

\caption{The flow chart of the BCG method\label{BCGflow}
}
\end{figure}

The BCG method is presented as follows (see Figure \ref{BCGflow}). First, the involved parameters and sets presented in Table \ref{Notations} and the stop criterion parameters (upper bound $UB$, lower bound $LB$, and tolerance $\epsilon$) are initialized. Note that the initialization of $P^r$ and $P^{a,t}$ significantly affects the efficiency of the CG in solving the BMP and BSP. Some literature (\citeauthor{maher2016solving}, \citeyear{maher2016solving}; \citeauthor{liang2018column}, \citeyear{liang2018column}) generates the initial aircraft routes according to the planned routes and the shortest paths. In this paper, we use constraint programming (CP) to generate aircraft routes and gate patterns since the CP can easily formulate the logic relations embedded in routes and patterns. In addition, commercial CP solvers (e.g., CPLEX) can solve CP models efficiently (\citeauthor{benoist2002constraint}, \citeyear{benoist2002constraint}). For more information about CP, we refer interested readers to \citeauthor{apt2003} (\citeyear{apt2003}) and \citeauthor{rossi2006} (\citeyear{rossi2006}). The CP models generating aircraft routes and gate patterns are provided in \ref{appendix CP models}.\par

Next, as shown in \textit{\textbf{RxBMP Optimality}} in Figure \ref{BCGflow}, we deal with the optimality of the RxBMP by Algorithm \ref{alg1}. Note that $(y^{r}_{p}, z_{j}, q) = (\mathbf{0}, \mathbf{1}, 0)$ is a feasible solution for the restricted problem of RxBMP (RtRxBMP) during the iteration of the BCG method, where $\mathbf{0}$ and $\mathbf{1}$ denote the vectors whose components are all equal to zero and one respectively. Hence, the RtRxBMP and its linear relaxation LR-RtRxBMP are always feasible. Once an aircraft route with negative reduced cost does not exist, the lower bound $LB$ is updated by the optimal value of the current LR-RtRxBMP, denoted by $v_{LM}^{*}$. Then, we check whether the optimality gap $(UB-LB)/LB\times100\%$ satisfies a given tolerance $\epsilon$.\par

\begin{algorithm}
\small
\setlength{\baselineskip}{1\baselineskip}
\caption{
Column Generation Algorithm for the Optimality of the RxBMP.
}
\label{alg1}
\flushleft{\vspace{-0.85em}{\bf Input:} the sets $\bar{\Pi}_{Fea}$ and $\bar{\Pi}_{Opt}$.}
\flushleft{\vspace{-0.8em}{\bf Output:} a flight and aircraft recovery decision $(y_{p}^{r*}, z_{j}^{*})$ and the updated values $(q^{*}, v_{M}^{*}, LB)$.}\\
\begin{algorithmic}[1]
\STATE Solve the LR-RtRxBMP and update its optimal value $v_{LM}^{*}$.
\WHILE{there exists an aircraft route with negative reduced cost}
    \FOR{each aircraft $r \in R$}
        \IF{there exists an aircraft route $p \in P^r$ whose reduced cost $\bar{c}_p^{r}$ is negative}
            \STATE Add the route $p$ to the set $\bar{P}^{r}$.
        \ENDIF
    \ENDFOR
    \STATE Solve the LR-RtRxBMP and update its optimal value $v_{LM}^{*}$.
\ENDWHILE
\STATE Solve the RtRxBMP and obtain the flight and aircraft recovery decision $(y_{p}^{r*}, z_{j}^{*})$.
\STATE Update the optimal value of the RtRxBMP $v_{M}^{*}$, the value $q^{*}$, and the lower bound $LB = \max \{LB, v_{LM}^{*}\}$.
\RETURN the flight and aircraft recovery decision $(y_{p}^{r*}, z_{j}^{*})$ and the values $(q^{*}, v_{M}^{*}, LB)$.
\end{algorithmic}
\end{algorithm}

Given an RxBMP solution $y_{p}^{r*}$, the feasibility of the associated BSP is checked first by Algorithm \ref{alg2}, denoted by \textit{\textbf{BSP Feasibility}} in Figure \ref{BCGflow}. Recall that the BSP is solved by CG, so its feasibility is not equivalent to the feasibility of the restricted problem of the BSP (RtBSP). When the RtBSP is infeasible, the \textit{phase one optimization problem} of its linear relaxation (P1LR-RtBSP) is introduced to deal with the feasibility of the BSP, which is well-known to tackle the feasibility in optimization problems (\citeauthor{boyd2004convex}, \citeyear{boyd2004convex}). When the LR-RtBSP is feasible but the Rt-BSP is still infeasible, we reject $ y_{p}^{r*}$ and resolve the RtRxBMP to obtain other solutions (\citeauthor{petersen2012optimization}, \citeyear{petersen2012optimization}).\par

\begin{algorithm}
\small
\setlength{\baselineskip}{1\baselineskip}
\caption{
Column Generation Algorithm for the Feasibility of the BSP.
}
\label{alg2}
\flushleft{\vspace{-0.8em}{\bf Input:} the RxBMP solution $(y_{p}^{r*}, z_{j}^{*})$.}
\flushleft{\vspace{-0.8em}{\bf Output:} the feasibility of the BSP related to the RxBMP solution $(y_{p}^{r*}, z_{j}^{*})$.}\\
\begin{algorithmic}[1]
\IF{the RtBSP related to $y_{p}^{r*}$ is feasible}
    \RETURN the fact that the BSP is feasible.
\ENDIF
\STATE Solve the P1LR-RtBSP (phase one optimization problem of the linear relaxation of the BSP).
\WHILE{there exists a gate pattern with negative reduced cost}
    \FOR{each gate type $t \in T_{a}$ at airport $a \in A$}
        \IF{there exists a gate pattern $p \in P^{a,t}$ whose reduced cost $\bar{c}_p^{a,t}$ is negative}
            \STATE Add the pattern $p$ to the set $\bar{P}^{a,t}$.
        \ENDIF
    \ENDFOR
    \IF{the RtBSP related to $y_{p}^{r*}$ is feasible}
        \RETURN the fact that the BSP is feasible.
    \ENDIF
    \STATE Solve the P1LR-RtBSP.
\ENDWHILE
\IF{the LR-RtBSP is feasible}
    \STATE Reject the solution $(y_{p}^{r*}, z_{j}^{*})$ in the RxBMP.
\ELSE
    \STATE Obtain an extreme direction $\pi_{d}$ from the dual feasible region of the P1LR-RtBSP.
    \STATE Update the set $\bar{\Pi}_{Fea} = \bar{\Pi}_{Fea} \bigcup \{\pi_{d}\}$.
\ENDIF
\RETURN the BSP is infeasible.
\end{algorithmic}
\end{algorithm}

If the BSP related to $y_{p}^{r*}$ is feasible, we move on to check its optimality by  Algorithm \ref{alg3}, \textit{\textbf{BSP Optimality}} as shown in Figure \ref{BCGflow}. The BCG method tackles optimality by applying CG to solve the linear relaxation of the BSP (LR-BSP). Once a gate pattern with negative reduced cost does not exist, the upper bound $UB$ is updated by the value $\left( v_{M}^{*} - q^{*} + v_{S}^{*} \right)$, where $v_{S}^{*}$ denotes the optimal value of the RtBSP. Then, we check whether the optimality gap is less than or equal to a given tolerance. If the BCG method is stopped, and the airline's integrated recovery decision is obtained from the solution corresponding to the current $UB$.\par

\begin{algorithm}
\small
\setlength{\baselineskip}{1\baselineskip}
\caption{
Column Generation Algorithm for the Optimality of the BSP.
}
\label{alg3}
\flushleft{\vspace{-0.8em}{\bf Input:} the RxBMP solution $(y_{p}^{r*}, z_{j}^{*})$, the values $(q^{*}, v_{M}^{*}, LB)$, and the tolerance $\epsilon$.}
\flushleft{\vspace{-0.8em}{\bf Output:} the updated upper bound $UB$.}\\
\begin{algorithmic}[1]
\STATE Solve the LR-RtBSP.
\WHILE{there exists a gate pattern with negative reduced cost}
    \FOR{each gate type $t \in T_{a}$ at airport $a \in A$}
        \IF{there exists a gate pattern $p \in P^{a,t}$ whose reduced cost $\bar{c}_p^{a,t}$ is negative}
            \STATE Add the pattern $p$ to the set $\bar{P}^{a,t}$.
        \ENDIF
    \ENDFOR
    \STATE Solve the LR-RtBSP.
\ENDWHILE
\STATE Solve the RtBSP and update its optimal value $v_{S}^{*}$ and the upper bound $UB = \max \{UB, v_{M}^{*}-q^{*}+v_{S}^{*}\}$.
\IF{$(UB-LB)/LB > \epsilon$}
    \STATE Obtain an extreme point $\pi_{p}$ from the dual feasible region of the LR-RtBSP.
    \STATE Update the set $\bar{\Pi}_{Opt} = \bar{\Pi}_{Opt} \bigcup \{\pi_{p}\}$.
\ENDIF
\RETURN the updated upper bound $UB$.
\end{algorithmic}
\end{algorithm}

Since the BSP is solved by CG, its feasibility and optimality usually cannot be detected before executing a complete process of CG. In addition, a complete process of CG requires considerable time for a large-scale problem. \citeauthor{petersen2012optimization} (\citeyear{petersen2012optimization}) developed two certificates (Theorem 4.1 and 4.2 in their work) to deal with this challenge by judging the infeasibility and suboptimality of their BSP during the CG process. Compared with the BSP in \citeauthor{petersen2012optimization} (\citeyear{petersen2012optimization}), which includes only linear equality constraints, the BSP in our work further considers the capacity inequality constraints (\ref{GRMcon4}). Inspired by  \citeauthor{petersen2012optimization} (\citeyear{petersen2012optimization}), we propose a certificate in Section \ref{Infeasibility Certificate} to judge the infeasibility of our BSP, which can enhance the efficiency of the BCG method.\par

Some literature (e.g., \citeauthor{barnhart1998flight}, \citeyear{barnhart1998flight}; \citeauthor{maher2016solving}, \citeyear{maher2016solving}) focused on developing column generation-based exact algorithms (e.g., branch-and-price) to solve mixed integer linear programming (MILP). Exact algorithms usually require a lot of computational costs to get an optimal solution, especially for large-scale MILPs. Due to the real-time requirement of solving airline integrated recovery problems, inspired by \citeauthor{petersen2012optimization} (\citeyear{petersen2012optimization}) and \citeauthor{liang2018column} (\citeyear{liang2018column}), we apply the CG strategy that solves the linear relaxation of MILPs by CG to choose candidate columns and then obtains an integer solution by solving the MILP with these candidate columns. Since the lower bound $LB$ and upper bound $UB$ in our BCG method are updated by the optimal value of the linear relaxation of RxBMP (the relaxation of SAGRM) and a feasible solution of SAGRM respectively, the above CG strategy would not affect the validity of $LB$ and $UB$.\par

\section{Acceleration techniques for the BCG method}\label{Section4}
In real-world airline disruption management, the SAGRM is required to be solved efficiently. The BMP, an aircraft rerouting problem, has been well studied in the existing literature (\citeauthor{liang2018column}, \citeyear{liang2018column}). However, the BSP, a gate reassignment prolem in an airline network, has not been studied in the airline industry and is often a large-scale problem. It means that the efficiency of solving the BSP plays an important role in dealing with the airline integrated recovery with gate reassignment. Hence, the following two acceleration techniques based on the particular characteristics and structure of the BSP are developed to improve the efficiency of the proposed BCG method.\par

\subsection{Separation technique}\label{Separation}
Separation is a favorable characteristic for large-scale optimization problems. If one large-scale problem has a separable structure, parallel computing or other advanced techniques can be applied to reduce the computational cost. It is noted that the BSP (\ref{GRMobj})$-$(\ref{GRMcon5}), a gate reassignment problem for all airports in an airline network, is not separable. However, under the Benders decomposition framework, the following proposition shows that the BSP can be separated into several small-scale gate reassignment problems related to each available gate type at each airport. Before providing the proposition, we give a lemma first to obtain the potential separable structure in the BSP, which could be used to separate the BSP.\par

\begin{lemma}\label{Lemma 1}
The set of flight copies $FC$ can be divided into several disjoint subsets of $FC$ according to the attributes of flight copies (the associated original flight, the departure/arrival airport, and the original fleet type) as follows:

\vspace*{-0.25cm}
\begin{align}
\label{lemma1_F}
FC = \bigcup_{j \in F} FC_j = \bigcup_{a \in A, t \in T_{a}} FC^{a,t}_{dep} = \bigcup_{a \in A, t \in T_{a}} FC^{a,t}_{arr}.
\end{align}

In addition, the set of flight copy connections $E$ can be divided into several disjoint subsets of $E$ according to the attributes of the connections as follows:

\vspace*{-0.25cm}
\begin{align}
\label{lemma1_E}
E = \bigcup\limits_{a \in A, t \in T_{a}} E^{a, t}.
\end{align}
\end{lemma}
Based on Lemma \ref{Lemma 1}, we give the following proposition to separate the BSP.\par

\begin{proposition}\label{Proposition 3}
Given an RxBMP solution ${\bar{y}_p^r}$, the BSP can be separated into $\sum\limits_{a \in A}|T_{a}|$ subproblems related to each available gate type at each airport. The separated subproblem for airport $a \in A$ and gate type $t \in T_{a}$ is
\begin{eqnarray}
&\min &\sum_{p\in {P^{a,t}}}c_{a,t} w_p^{a,t} \label{SPobjective}\\
&\mathrm{s.t.}&\sum\limits_{p\in {P^{a,t}}}a_{{p,{j_{v}dep}}}^{a,t} w_p^{a,t}
     = \sum\limits_{r\in R}\sum\limits_{p\in {P^r}}a_{{{p},j_v}}^{r} \bar{y}_p^r,\ \ \forall {j_v\in {FC^{a,t}_{dep}}}, \label{SPConstraint1}\\
&&
\sum\limits_{p\in {P^{a,t}}}a_{{p,{j_{v}arr}}}^{a,t} w_p^{a,t}
     = \sum\limits_{r\in R}\sum\limits_{p\in {P^r}}a_{{{p},j_v}}^{r} \bar{y}_p^r,\ \ \forall {j_v\in {FC^{a,t}_{arr}}}, \label{SPConstraint2}\\
&&
\sum\limits_{p\in {P^{a,t}}}b_{{p,{(i_u,j_v)}}}^{a,t} w_p^{a,t} = 1,\ \ \forall (i_u,j_v)\in {E^{a, t}}\cap E_1, \label{SPConstraint3}\\
&&
\sum\limits_{p\in {P^{a,t}}} w_p^{a,t} \leq n_{a,t},  \label{SPConstraint4}\\
&&
w_p^{a,t}\in \{0,1\},\ \ \forall  p\in P^{a,t},  \label{SPConstraint5}
\end{eqnarray}
where $E_1$ denotes the set $\left\{(i_u,j_v)~{\Bigg|}~\sum\limits_{r\in R}\sum\limits_{p\in {P^r}}b_{p,(i_u,j_v)}^{r} \bar{y}_p^r = 1\right\}$.
\end{proposition}

\begin{proof}
With the aid of Equation (\ref{lemma1_F}), Constraints (\ref{GRMcon1}) and (\ref{GRMcon2}) can be divided according to the departure/arrival airports of flight copies and the associated gate types as follows:

\vspace*{-0.5cm}
\begin{align}
    \nonumber
  &\sum_{a\in A}\sum_{t\in T_{a}}\sum_{p\in {P^{a,t}}}a_{{p,{j_{v}dep}}}^{a,t} w_p^{a,t}
     = \sum_{r\in R}\sum_{p\in {P^r}}a_{{{p},j_v}}^{r} \bar{y}_p^r,\ \ \forall j\in {F},{j_v\in {FC_j}}\\
     \iff
  &\sum_{p\in {P^{a,t}}}a_{{p,{j_{v}dep}}}^{a,t} w_p^{a,t}
     = \sum_{r\in R}\sum_{p\in {P^r}}a_{{{p},j_v}}^{r} \bar{y}_p^r,\ \ \forall {j_v\in {FC^{a,t}_{dep}}}, {a\in A}, {t\in T_{a}}, \label{Separation1}
\end{align}
and

\vspace*{-0.5cm}
\begin{align}
  \nonumber
  &\sum_{a\in A}\sum_{t\in T_{a}}\sum_{p\in {P^{a,t}}}a_{{p,{j_{v}arr}}}^{a,t} w_p^{a,t}
     = \sum_{r\in R}\sum_{p\in {P^r}}a_{{{p},j_v}}^{r} \bar{y}_p^r,\ \ \forall j\in {F},{j_v\in {FC_j}}\\
     \iff
  &\sum_{p\in {P^{a,t}}}a_{{p,{j_{v}arr}}}^{a,t} w_p^{a,t}
     = \sum_{r\in R}\sum_{p\in {P^r}}a_{{{p},j_v}}^{r} \bar{y}_p^r,\ \ \forall {j_v\in {FC^{a,t}_{arr}}}, {a\in A}, {t\in T_{a}}. \label{Separation2}
\end{align}
Due to Equation (\ref{lemma1_E}), Constraints (\ref{GRMcon3}) can be divided according to the connecting airport and gate type, i.e.,

\vspace*{-0.5cm}
\begin{align}
  \nonumber
  &\sum_{a\in A}\sum_{t\in T_{a}}\sum_{p\in {P^{a,t}}}b_{{p,{(i_u,j_v)}}}^{a,t} w_p^{a,t}
     \geq \sum_{r\in R}\sum_{p\in {P^r}}b_{p,(i_u,j_v)}^{r} \bar{y}_p^r ,\ \ \forall (i_u,j_v)\in {E}\\
    \overset{\eqref{lemma1_E}}\iff
    \nonumber
  &\sum_{p\in {P^{a,t}}}b_{{p,{(i_u,j_v)}}}^{a,t} w_p^{a,t}
     \geq \sum_{r\in R}\sum_{p\in {P^r}}b_{p,(i_u,j_v)}^{r} \bar{y}_p^r ,\ \ \forall (i_u,j_v)\in {E^{a, t}}, {a\in A}, {t\in T_{a}} \\
    \iff
  \nonumber
  &\sum_{p\in {P^{a,t}}}b_{{p,{(i_u,j_v)}}}^{a,t} w_p^{a,t}
     \geq \sum_{r\in R}\sum_{p\in {P^r}}b_{p,(i_u,j_v)}^{r} \bar{y}_p^r ,\ \ \forall (i_u,j_v)\in {E^{a, t}}\cap E_1, {a\in A}, {t\in T_{a}},
\end{align}
where $E_1$ denotes the set $\{(i_u,j_v)| \sum\limits_{r\in R}\sum\limits_{p\in {P^r}}b_{p,(i_u,j_v)}^{r} \bar{y}_p^r = 1\}$. Since $\sum\limits_{a\in A}\sum\limits_{t\in T_{a}}\sum\limits_{p\in {P^{a,t}}}b_{{p,{(i_u,j_v)}}}^{a,t} w_p^{a,t}$ is binary for any $(i_u,j_v) \in E$, we obtain that

\vspace*{-0.5cm}
\begin{align}
  \nonumber
  &\sum_{p\in {P^{a,t}}}b_{{p,{(i_u,j_v)}}}^{a,t} w_p^{a,t}
     \geq \sum_{r\in R}\sum_{p\in {P^r}}b_{p,(i_u,j_v)}^{r} \bar{y}_p^r ,\ \ \forall (i_u,j_v)\in {E^{a, t}}\cap E_1, {a\in A}, {t\in T_{a}}\\
  \iff
  &\sum_{p\in {P^{a,t}}}b_{{p,{(i_u,j_v)}}}^{a,t} w_p^{a,t} = 1,\ \ \forall (i_u,j_v)\in {E^{a, t}}\cap E_1, {a\in A}, {t\in T_{a}}.
     \label{Separation3}
\end{align}

Note that the objective function (\ref{GRMobj}) and the remaining Constraints (\ref{GRMcon4})$-$(\ref{GRMcon5}) are also separable in airport $a$ with gate type $t$. Combining (\ref{GRMobj}), (\ref{GRMcon4})$-$(\ref{GRMcon5}) and (\ref{Separation1})$-$(\ref{Separation3}), the BSP can be separated into $\sum\limits_{a \in A}|T_{a}|$ subproblems related to airport $a \in A$ and gate type $t \in T_{a}$.
\end{proof}

\vspace{0.5em}
This proposition shows the separation characteristic in the BSP with the aid of the reformulation of Constraints (\ref{GRMcon1})$-$(\ref{GRMcon3}). By dividing the BSP into several subproblems related to each gate type at each airport, this separation technique greatly alleviates the difficulty of solving the BSP and then accelerates the efficiency of the BCG method.\par

\subsection{Infeasibility certificate}\label{Infeasibility Certificate}
Note that not all columns are involved in the LR-RtBSP. Then, if the LR-RtBSP is infeasible, the feasibility of the LR-BSP cannot be determined unless a complete CG process is carried out or a feasible solution is found during the CG process. However, executing a complete CG will take a lot of efforts if the size of the BSP is large. \citeauthor{petersen2012optimization} (\citeyear{petersen2012optimization}) proposed an infeasibility certificate to consider the linear equality constraints in their Benders subproblem when it was solved by CG . Motivated by their infeasibility certificate, we propose an infeasibility certificate (Theorem \ref{Theorem1}) that further considers the capacity inequality constraints (\ref{GRMcon4}) in the Benders subproblem. Once the condition of the theorem is satisfied, we can stop the CG of LR-BSP and add a no-good cut to RxBMP. The following proposition is built first before providing Theorem \ref{Theorem1}.\par

\begin{proposition}
\label{Proposition 4}
Let $A\in \mathbf{R}^{m \times n}, b\in \mathbf{R}^{m}$ and $t\in \mathbf{R}$. One and only one of the following two sets

\vspace*{-0.5cm}
$$ X:=\{ x~|~Ax = b, \mathbf{1}^{T}x \leq t ,  x\succeq0 \} $$
and

\vspace*{-0.5cm}
$$\Lambda_{\delta}:=\{ (\lambda,\delta)~|~\lambda^{T}A-\delta \mathbf{1}^{T} \preceq 0, \lambda^{T}b-t\delta>0, \delta\geq0 \}$$
is nonempty.
\end{proposition}

\begin{proof}
By introducing a slack variable $x_{t}\geq0$, the set $X$ is transformed into the following standard form

\vspace*{-0.5cm}
$$X_{t}:=\{ (x,x_t)~|~Ax = b, -(\mathbf{1}^{T}x + x_t) = -t,  (x,x_t)\succeq0 \}.$$
It is clear that $X_t$ is nonempty if and only if $X$ is nonempty.\par

By Farkas' lemma (\citeauthor{boyd2004convex}, \citeyear{boyd2004convex}), it is known that the linear systems corresponding to $X_{t}$ and $\Lambda_{\delta}$ are strong alternatives, which leads to the result that only one of the sets $X$ and $\Lambda_{\delta}$ is nonempty.\end{proof}\par

Based on Proposition \ref{Proposition 4}, we present the next proposition which will be used to prove Theorem \ref{Theorem1}.\par

\begin{proposition}
\label{Proposition 5}
Let $\bar{A}$ be a matrix that satisfies $col(\bar{A})\subseteq col(A)$, where $col(\cdot)$ denotes the set of the columns of a given matrix. Suppose there exists $(\bar{\lambda},\bar{\delta})$  with $\bar{\delta}\geq0$ satisfying  $\bar{\lambda}^{T}a-\bar{\delta} \leq 0$ for any vector $a \in col(\bar{A})$.
If $(\bar{\lambda},\bar{\delta})$ also satisfies $\bar{\lambda}^{T}b-t\bar{\delta}>0$ and $\bar{\lambda}^{T}a-\bar{\delta} \leq 0$ for any vector $a \in col(A) \backslash col(\bar{A})$, then $X$ is empty.
\end{proposition}

\begin{proof}
It is clear that $(\bar{\lambda},\bar{\delta}) \in \Lambda_{\delta}$, which means $\Lambda_{\delta}$ is nonempty.
Then, it follows from Proposition \ref{Proposition 4} that $X$ is empty.\
\end{proof}

After applying the separation technique in Section \ref{Separation}, the feasible region of the linear relaxation of the separated subproblem related to airport $a$ and gate type $t$ (see Equtions \eqref{SPConstraint1}$-$\eqref{SPConstraint5}) can be transformed into the following set

\vspace*{-0.4cm}
$$X_{a,t}:=\{ x~|~Ax = b_{\bar{y}_{p}^{r}}, \mathbf{1}^{T}x \leq n_{a,t} ,  x\succeq0 \},$$
where $Ax = b_{\bar{y}_{p}^{r}}$ and $\mathbf{1}^{T}x \leq n_{a,t}$ represent Constraints \eqref{SPConstraint1}$-$\eqref{SPConstraint3} and the gate capacity constraint \eqref{SPConstraint4}, respectively. Note that each column of $A$ corresponds to one possible gate pattern $p \in P^{a,t}$. During the iterative process of CG, only the part of patterns $\bar{P}^{a,t} \subseteq P^{a,t}$  in the LR-RtBSP, and $\bar{A}$ is set to be the matrix composed by these columns. Then $X_{a,t}$ is restricted to the set

\vspace*{-0.4cm}
$$\bar{X}_{a,t}:=\{ \bar{x}~|~\bar{A}\bar{x} = b_{\bar{y}_{p}^{r}}, \mathbf{1}^{T}\bar{x} \leq n_{a,t},  \bar{x}\succeq0 \}.$$
According to Proposition \ref{Proposition 5}, the following theorem gives a condition to judge whether $X_{a,t}$ is empty when $\bar{X}_{a,t}$ is empty.\par

\begin{theorem}
\label{Theorem1}
Suppose $\bar{X}_{a,t}$ is empty for the given RxBMP solution $\bar{y}_{p}^{r}$.
Let $\bar{\lambda}$ and $\bar{\delta}$ denote the dual variables corresponding to $\bar{A}\bar{x}=b_{\bar{y}_{p}^{r}}$ and $\mathbf{1}^{T}\bar{x} \leq n_{a,t}$ in the P1LR-RtBSP respectively, and let $\hat{\delta} = \textrm{max}\left\{ \bar{\delta},\textrm{max}\{\bar{\lambda}^{T}a|a \in col{(A)}\setminus col(\bar{A})\}\right\}$.
If $\bar{\lambda}^{T}b_{\bar{y}_{p}^{r}}-n_{a,t}\hat{\delta}>0$, then $X_{a,t}$ is empty and the BSP related to $\bar{y}_{p}^{r}$ is infeasible.
\end{theorem}

\begin{proof}
We only need to prove  $(\bar{\lambda},\hat{\delta})$ satisfies the conditions of Proposition \ref{Proposition 5}, i.e.,
\begin{equation} \label{setrelation}
(\bar{\lambda},\hat{\delta}) \in \{ (\lambda,\delta)~|~\lambda^{T}A-\delta \mathbf{1}^{T} \preceq 0, \lambda^{T}b_{\bar{y}_{p}^{r}}-n_{a,t}\delta>0, \delta\geq0 \}.
\end{equation}

The last two inequalities, i.e., $~\bar{\lambda}^{T}b_{\bar{y}_{p}^{r}}-n_{a,t}\hat{\delta}>0$ and
$\hat{\delta}\geq 0$,  holds immediately from the condition of the theorem.

Now we prove the first inequality, i.e., $\bar{\lambda}^{T}a-\hat{\delta} \leq 0$ holds for any column $a \in col{(A)}$. Note that the reduced cost corresponding to column $a$ is $-\bar{\lambda}^{T}a+\bar{\delta}$.
    \begin{itemize}
    \vspace*{0cm}
    \item For any column $a \in col(\bar{A})$, $\bar{\lambda}^{T}a-\hat{\delta} \leq \bar{\lambda}^{T}a-\bar{\delta} \leq 0$ holds since $\hat{\delta} \geq \bar{\delta}$.
    \vspace*{-0.1cm}
    \item For any column $a \in col{(A)}\setminus col(\bar{A})$, $\bar{\lambda}^{T}a-\hat{\delta} \leq 0$ holds because
    \vspace*{-0.1cm}
            \begin{center}$\hat{\delta} \geq \textrm{max}\{\bar{\lambda}^{T}a~|~a \in col{(A)}\setminus col(\bar{A})\}\}.$\end{center}
    \end{itemize}

Summing up, the conditions of Proposition \ref{Proposition 5} hold, then it follows that $X_{a,t}$ is empty.
Hence, the BSP related to $\bar{y}_{p}^{r}$ is infeasible since the linear system corresponding to $X_{a,t}$ (the subset of all constraints in the BSP) is infeasible.
\end{proof}

Note that Theorem \ref{Theorem1} is a sufficient condition to judge the infeasibility. It means that $\bar{y}_{p}^{r}$ might be infeasible for the BSP (with all possible patterns in $P^{a,t}$) when the condition in the theorem does not hold. In this case, the CG is required to be carried on until we can judge the infeasibility by the theorem or the CG is completed. However, Theorem \ref{Theorem1} provides a chance to terminate the CG process and promotes the efficiency of the BCG method.\par

\begin{remark}
With the aid of the separation technique in Section \ref{Separation}, detecting the infeasibility of the Benders subproblem \eqref{GRMobj}$-$\eqref{GRMcon5} is equivalent to detecting the infeasibility related to the separated problems of the Benders subproblem. Hence, once Theorem \ref{Theorem1} detects the infeasibility for one separated subproblem, the current CG process can be terminated immediately, and a no-good cut (\citeauthor{wolsey2020integer}, \citeyear{wolsey2020integer}) can be added to RxBMP.\par
\end{remark}

\section{Computational experiments and discussions}\label{Section5}
To validate the effectiveness of the proposed SAGRM and the BCG method, we carry out different kinds of experiments on the open dataset initially provided by the ROADEF 2009 challenge (\url{https://www.roadef.org/challenge/2009/en/instances.php}). All experiments are coded in Java and run on a laptop with a 2.60 GHz Intel i7 CPU and 12 threads. All models involved in the BCG method (LPs, MILPs, and CPs) are solved by CPLEX 12.9.\par

\subsection{Data description and experiment setting}
\noindent
As summarized in Table \ref{Characteristics of test instances}, ten instances with different disruption scales are extracted from the ROADEF 2009 challenge. These instances have two types of disruptions, namely hub closure and flow control, with varying lengths of disruptions. The recovery time window, similar to the case in \citeauthor{petersen2012optimization} (\citeyear{petersen2012optimization}), is set to be one day in the experiments. Besides, we consider two types of gates for wide and narrow-body planes. To simplify our discussions, we assume that narrow-body and wide-body aircraft are only allowed to be assigned to narrow and wide-body gates respectively, which can be easily extended to the scenario where wide-body gates are compatible with narrow-body aircraft.\par

\begin{table}[!htbp]
  \centering
  \renewcommand{\baselinestretch}{1.1}
  \caption{Characteristics of data set in computational experiments\label{Characteristics of test instances}}
  \normalsize
    {\resizebox{1\linewidth}{!}{
    \begin{tabular}{cccccc}
    \toprule
    Instance & No. of flights & No. of aircraft & No. of fleet types & No. of airports & Disruption \\
    \midrule
    1     & 24    & 4     & 2     & 7     & Hub closure 1h \\
    2     & 52    & 8     & 3     & 9     & Flow control 1h \\
    3     & 98    & 15    & 3     & 14    & Flow control 2h \\
    4     & 150   & 24    & 4     & 18    & Flow control 1h \\
    5     & 198   & 32    & 4     & 17    & Flow control 1h \\
    6     & 250   & 41    & 4     & 22    & Flow control 2h \\
    7     & 295   & 49    & 5     & 26    & Flow control 2h \\
    8     & 359   & 61    & 6     & 28    & Flow control 2h \\
    9     & 395   & 68    & 8     & 30    & Hub closure 1h \& Flow control 1h \\
    10    & 464   & 81    & 11    & 35    & Hub closure 1h \& Flow control 1h \\
    \bottomrule
    \end{tabular}}}
\end{table}

The original gate assignment in the ROADEF 2009 challenge is not provided. Thus, we generated the required gate capacity parameter from the original flight schedules as follows. First, we obtain a gate assignment solution corresponding to the minimum total number of required gates by solving the GRM (\ref{GRMobj})$-$(\ref{GRMcon5}) without Constraints (\ref{GRMcon4}), where $\bar{y}_p^r$ is set according to the original flight schedule. Let $\tilde{n}_{a,t}$ denote the number of used gates at airport $a\in A$ with gate type $t\in T_{a}$ in the above problem. Then, the associated number of gates $n_{a, t}$ is estimated as $n_{a, t}  = \lceil (1+10\%)\cdot \tilde{n}_{a,t} \rceil$, where the 10\% correction factor accounts for the gates that are available but not used in the original plan.\par

Key parameters used in our experiments are given in Table \ref{Benchmark parameters used in experiments}, which are based on existing literature. All delay intervals in Table \ref{Benchmark parameters used in experiments} are less than that in \citeauthor{maher2016solving} (\citeyear{maher2016solving}), and we expect the recovery plans obtained under this condition to be economical. The effect of adopting different delay intervals on the runtime and solution quality will be further discussed in Section \ref{Section The length of delay interval and gate types}. The cost parameters for flight cancellation, delay, and swap are adopted from \citeauthor{petersen2012optimization} (\citeyear{petersen2012optimization}). The cost of one gate pattern $c_{a,t}$ is estimated according to \citeauthor{zhang2017optimization} (\citeyear{zhang2017optimization}). In addition, the optimality gap for the BCG method, determined by $(UB-LB)/LB\times100\%$, is set as 5\%.\par

\begin{table}[!htbp]
  \centering
  \renewcommand{\baselinestretch}{1}
  \caption{Key parameters used in computational experiments\label{Benchmark parameters used in experiments}}
  \small
    \begin{tabular}{l|cccccccccc}
    \toprule
    \multicolumn{1}{l}{Instance} & 1     & 2     & 3     & 4     & 5     & 6     & 7     & 8     & 9     & 10 \\
    \midrule
    \cmidrule{2-11} Length of delay interval (mins) & 5     & 5     & 5     & 5     & 5     & 5     & 10    & 10    & 15    & 15 \\
    Cost of cancellation per flight (unit\footnotemark[1]) & \multicolumn{10}{c}{200} \\
    Cost of flight delay per minute (unit) & \multicolumn{10}{c}{1} \\
    Cost of gate usage (unit) & \multicolumn{10}{c}{4} \\
    Gate types & \multicolumn{10}{c}{2} \\
    Maximum allowable delay (mins) & \multicolumn{10}{c}{120} \\
    \bottomrule
    \end{tabular}\\
     \raggedright
     {\hspace*{1.5cm} \footnotemark[1]{\footnotesize 1 unit = USD\$125}}
\end{table}

\subsection{Computational results}\label{Section Computational Results}
\noindent
In order to evaluate the performance of the BCG method, test instances are also solved by the well-known sequential recovery approach (SEQ), i.e.,  sequentially solving the SARM (\ref{SARMobj})$-$(\ref{SARMcon9}) without Benders cuts (\ref{classFeaCut})$-$(\ref{classOptCut}) and the GRM (\ref{GRMobj})$-$(\ref{GRMcon5}). To ensure a fair comparison, the separation technique in Section \ref{Separation} are also applied to the SEQ method. \par

As mentioned in the first motivation of considering the integrated recovery problem in Section \ref{Introduction}, the SEQ method considers the gate capacity restriction (\ref{GRMcon4}) by reducing the slot capacity $u_{sarr}$ in Constraints (\ref{SARMcon3}) for the schedule and aircraft recovery problem when the gate resource is not sufficient. Let $\alpha$ denote the rate of reduction of the slot capacity, and $u_{sarr}^{\textrm{SEQ}}$ denote the updated slot capacity that is calculated by

\vspace*{-0.5cm}
\begin{align}
u_{sarr}^{\textrm{SEQ}} =\lfloor (1-\alpha) \cdot u_{sarr}\rfloor.\label{u_alpha}
\end{align}
Recall that the estimation of the updated slot capacity in the SEQ method may lead to overestimation (SEQ-OE) and underestimation (SEQ-UE). To show the effect of these two estimations, the values of $\alpha$ in experiments are set within the ranges 3-10\% and 6-20\% in SEQ-OE and SEQ-UE, respectively. As shown in Table \ref{alpha value}, the specific value of $\alpha$ is set according to the scale of the corresponding instance.\par

\begin{table}[!htbp]
  \centering
  \renewcommand{\baselinestretch}{1}
  \caption{The value of $\alpha$ in \eqref{u_alpha} for the sequential recovery approach\label{alpha value}}
  \small
    \begin{tabular}{l|cccccccccc}
    \toprule
    \multicolumn{1}{l}{Instance} & 1     & 2     & 3     & 4     & 5     & 6     & 7     & 8     & 9     & 10 \\
    \midrule
    The value of $\alpha$ in \eqref{u_alpha} for SEQ-OE (\%) & 10     & 10     & 10     & 10     & 5     & 5     & 4    & 4    & 4    & 3 \\
    The value of $\alpha$ in \eqref{u_alpha} for SEQ-UE (\%) & 20     & 20     & 15     & 11     & 10     & 8     & 6    & 6    & 6    & 6 \\
    \bottomrule
    \end{tabular}
\end{table}

Table \ref{Results of integrated recovery and sequential recovery} provides the results of solving the test instances by the BCG method and the SEQ method. In addition, we also use CPLEX to solve the integrated recovery model SAGRM. However, it is found that CPLEX fails to solve these instances within reasonable computational times. For large-scale instances, it is even out of memory space when inputting the SAGRM formulation. Thus, we do not report the results of CPLEX.

\begin{table}[!htbp]
  \centering
  \renewcommand{\baselinestretch}{1.2}
  \caption{Results of sequential recovery SEQ and integrated recovery SAGRM\label{Results of integrated recovery and sequential recovery}}
  \small
    {\resizebox{\linewidth}{!}{
    \begin{tabular}{l|cccccccccc}
    \toprule
    Instance & 1     & 2     & 3     & 4     & 5     & 6     & 7     & 8     & 9     & 10 \\
    \midrule
    \multicolumn{11}{l}{Solution metrics by SEQ-OE} \\
    \midrule
    No. of canceled flights & 2     & 0     & 0     & 1     & 2     & 2     & 3     & 4     & 7     & 9 \\
    Total delay time (mins) & 20    & 95    & 172   & 100   & 95    & 250   & 420   & 460   & 510   & 435 \\
    No. of swapped tail assignments & 10    & 22    & 42    & 38    & 112   & 83    & 98    & 154   & 152   & 168 \\
    No. of used gates & 16    & 27    & {\bf infea} & 66    & 108   & 108   & 130   & 164   & {\bf infea} & {\bf infea} \\
    Schedule and aircraft recovery cost (unit) & 420   & 95    & 175   & 300   & 495   & 650   & 1020  & 1260  & 1910  & 2235 \\
    Gate recovery cost (unit) & 64    & 108   & {\bf infea} & 264   & 432   & 432   & 520   & 656   & {\bf infea} & {\bf infea} \\
    Total recovery cost (unit)  & 484   & 203   & {\bf infea} & 564   & 927   & 1082  & 1540  & 1916  & {\bf infea} & {\bf infea} \\
    CPU time (s) & 3     & 22    & 44    & 54    & 69    & 124   & 155   & 163   & 179   & 196 \\
    \midrule
    \multicolumn{11}{l}{Solution metrics by SEQ-UE} \\
    \midrule
    No. of canceled flights & 4     & 4     & 2     & 5     & 3     & 8     & 8     & 5     & 10    & 18 \\
    Total delay time (mins) & 0     & 65    & 150   & 20    & 95    & 70    & 150   & 450   & 510   & 285 \\
    No. of swapped tail assignments & 8     & 13    & 34    & 51    & 124   & 118   & 117   & 159   & 146   & 180 \\
    No. of used gates & 14    & 23    & 46    & 61    & 99    & 101   & 122   & 158   & 191   & 221 \\
    Schedule and aircraft recovery cost (unit) & 800   & 865   & 550   & 1020  & 695   & 1670  & 1750  & 1450  & 2510  & 3885 \\
    Gate recovery cost (unit) & 56    & 92    & 184   & 244   & 396   & 404   & 488   & 632   & 764   & 884 \\
    Total recovery cost (unit)  & 856   & 957   & 734   & 1264  & 1091  & 2074  & 2238  & 2082  & 3274  & 4769 \\
    CPU time (s) & 4     & 24    & 47    & 52    & 77    & 128   & 151   & 160   & 185   & 203 \\
    \midrule
    \multicolumn{11}{l}{Solution metrics by BCG} \\
    \midrule
    No. of canceled flights & 2     & 0     & 0     & 1     & 2     & 2     & 3     & 3     & 7     & 9 \\
    Total delay time (mins) & 20    & 95    & 195   & 100   & 100   & 250   & 420   & 660   & 510   & 435 \\
    No. of swapped tail assignments & 10    & 22    & 47    & 37    & 122   & 72    & 130   & 119   & 159   & 149 \\
    No. of used gates & 15    & 24    & 51    & 59    & 95    & 98    & 118   & 148   & 171   & 206  \\
    Schedule and aircraft recovery cost (unit) & 420   & 95    & 215   & 300   & 500   & 650   & 1020  & 1260  & 1910  & 2235 \\
    Gate recovery cost (unit) & 60    & 96    & 184   & 236   & 380   & 392   & 472   & 592   & 684   & 824 \\
    Total recovery cost (unit)  & 480   & 191   & 399   & 536   & 880   & 1042  & 1492  & 1852  & 2594  & 3059 \\
    Optimality gap (\%)  & 0.63\% & 0.70\% & 1.82\% & 1.01\% & 0.53\% & 0.80\% & 0.79\% & 0.84\% & 1.89\% & 0.52\% \\
    CPU time (s) & 4     & 34    & 62    & 76    & 118   & 192   & 258   & 269   & 282   & 294 \\
    \midrule
    \multicolumn{11}{l}{Savings by BCG comepared with SEQ-OE} \\
    \midrule
    Schedule and aircraft recovery saving (\%) & 0\%   & 0\%   & -23\% & 0\%   & -1\%  & 0\%   & 0\%   & 0\%   & 0\%   & 0\% \\
    Gate recovery saving (\%) & 6\%   & 11\%  & n/a   & 11\%  & 12\%  & 9\%   & 9\%   & 10\%  & n/a   & n/a \\
    Total recovery cost saving (\%) & 1\%   & 6\%   & n/a   & 5\%   & 5\%   & 4\%   & 3\%   & 3\%   & n/a   & n/a \\
    CPU time increase (\%) & 33\%  & 55\%  & 41\%  & 41\%  & 71\%  & 55\%  & 66\%  & 65\%  & 58\%  & 50\% \\
    \midrule
    \multicolumn{11}{l}{Savings by BCG comepared with SEQ-UE} \\
    \midrule
    Schedule and aircraft recovery saving (\%) & 48\%  & 89\%  & 61\%  & 71\%  & 28\%  & 61\%  & 42\%  & 13\%  & 24\%  & 42\% \\
    Gate recovery saving (\%) & -7\%  & -4\%  & 0\%   & 3\%   & 4\%   & 3\%   & 3\%   & 6\%   & 10\%  & 7\% \\
    Total recovery cost saving (\%) & 44\%  & 80\%  & 46\%  & 58\%  & 19\%  & 50\%  & 33\%  & 11\%  & 21\%  & 36\% \\
    CPU time increase (\%) & 0\%   & 42\%  & 32\%  & 46\%  & 53\%  & 50\%  & 71\%  & 68\%  & 52\%  & 45\% \\
    \bottomrule
    \end{tabular}}}
\end{table}

\subsection*{Insights from computational results}
We have the following observations from Table \ref{Results of integrated recovery and sequential recovery}. First, the computational times for solving test instances by the proposed BCG method are all less than 5 minutes with optimality gaps within 5\%. Compared with the SEQ methods, the runtime for the BCG method is longer for all test instances. The reason is that the integrated SAGRM model is larger and far more complex than sequentially solving the SARM and GRM. The advantage of the integrated SAGRM is that it avoids infeasible solutions encountered by the SEQ method while saving up to 80\% of total recovery costs in some scenarios. It is noted that infeasible gate reassignment issues result in further losses to airlines, airports, and passengers due to flight and aircraft gate delays.\par

Second, when the cost saving of the BCG method is broken down into savings from the schedule and aircraft recovery and the gate reassignment, it is seen in Table \ref{Results of integrated recovery and sequential recovery} that the schedule and aircraft recovery costs yielded by BCG are at least on par with those by SEQ-OE except Instances \#3 and \#5, while the cost saving from the gate reassignment by BCG is about 6-12\% than by SEQ-OE. This cost efficiency in the gate reassignment leads to lower total recovery costs for BCG in all test instances by 1-6\%. For the SEQ-UE scenario, due to the underestimation of the slot capacity. Hence, the schedule and aircraft recovery costs saving by the BCG method can be 13-89\% better than SEQ-UE. In addition, the cost savings from the gate reassignment by BCG is up to 10\% better than SEQ-UE, leading the overall cost saving between 11-80\% among test instances. \par

Third, the SEQ-OE method fails to generate a feasible solution for some instances, including  Instances \#3, \#9, and \#10. The main reason is that the SEQ-OE method overestimates the slot capacity, and this causes some arrival flights to be unable to be serviced by gates. The required total number of gates for Instances \#3, \#9, and \#10 by SEQ-OE are 64, 193, and 233, respectively if we relax the gate capacity restriction. This is still much higher than the gate usage by BCG, which is 51, 171, and 206 gates respectively for Instance \#3, \#9 and \#10. For Instance \#10, there are 27 gate shortfalls for the network. This large gate shortage will trigger a secondary gate assignment disruption and require significant effort to recover with extra operating costs. However, as shown in Table \ref{Results of integrated recovery and sequential recovery}, the BCG method can effectively avoid gate reassignment infeasibility, which demonstrates the superiority of our integrated recovery approach, especially when dealing with large-scale disruptions.\par

Fourth, SEQ-UE imposes higher restrictions on gate capacity to ensure the feasibility of gate reassignment. Although it avoids the infeasibility of gate reassignment by adopting a conservative slot estimation, this measure limits airline recovery options, causing longer flight delays. As shown in Table \ref{Results of integrated recovery and sequential recovery}, our SAGRM approach can still save from -7\% to 10\% gate recovery costs while maintaining huge savings in schedule and aircraft recovery. This results in total recovery cost savings from 11\% to 80\% among test instances. This is because the generated Benders cuts and valid inequalities in our BCG method  can capture the interdependence between the schedule, aircraft recovery and the gate recovery. This superiority helps the BCG method find a high-quality recovery plan through an integration perspective, which is beyond the ability of the sequential recovery approach.\par

Fifth, according to the results of Instances \#8, multiple feasible recovery plans might have the same recovery cost for the schedule and aircraft recovery, but have different effects on the gate recovery decisions. Specifically, Instance \#8 has the same schedule and aircraft recovery cost from the BCG method and the SEQ-OE method, 1260 units. However, the BCG method only requires 148 gates, while the SEQ-OE method needs 164 gates. Overall, the SEQ-OE method requires 7-14\% more gates to accommodate flights, and the SEQ-UE method still needs up to 12\% more gates than the integrated SAGRM approach.\par

\subsection{Computational techniques and their effects}
\subsubsection*{Separation technique}
\begin{figure}[!htbp]
\centering
\includegraphics[width=0.67\textwidth]{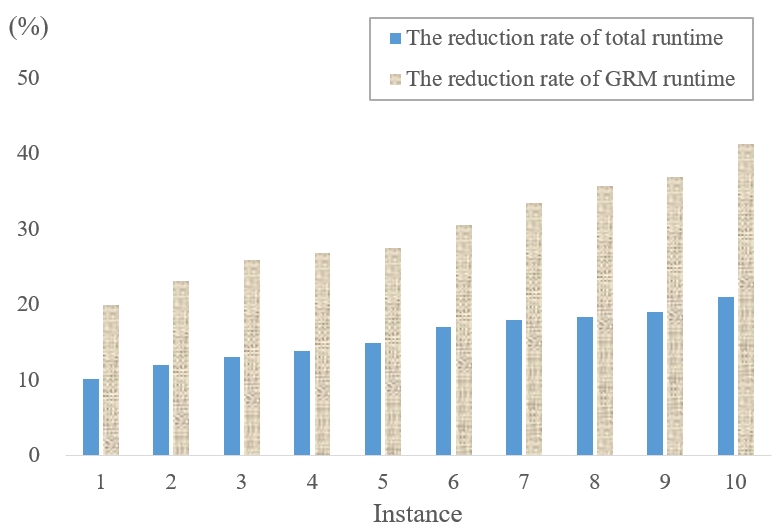}\\
\vspace{-0.0cm}
\caption{The reduction rate in runtime with separation technique\label{Sensitivity analysis for separation technique}}
\end{figure}

Figure \ref{Sensitivity analysis for separation technique} shows the effect of the separation technique proposed in Section \ref{Separation}. Compared to the scenario without this technique, the separation technique reduces the total runtime by 10-20\%. Such a reduction is attributed to the technique's ability to separate the Benders subproblem (GRM) into several small-scale gate reassignment problems, significantly lowering the complexity of solving the Benders subproblem. Then, with the aid of the parallel computing technique, the runtime is substantially reduced. The average and the maximum number of the reduction of GRM runtime are equal to 30\% and 41\% over these instances, respectively. In addition, it is clear that the above effect will be more pronounced as the number of airports and gate types increases. As shown in Figure \ref{Sensitivity analysis for separation technique}, the reduction rate becomes larger when the number of airports increases.\par

\subsubsection*{Infeasibility certificate}
\begin{figure}[!htbp]
\centering
\includegraphics[width=0.67\textwidth]{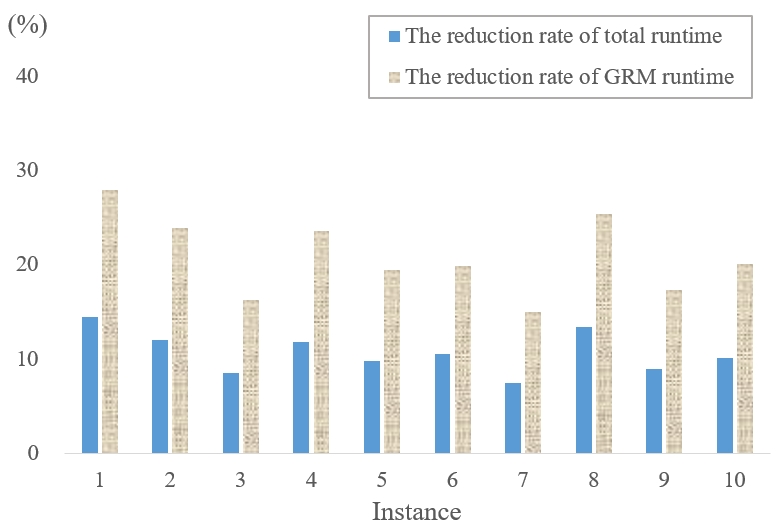}\\
\vspace{-0.0cm}
\caption{The reduction rate in runtime with infeasibility certificate \label{Sensitivity analysis for infeasibility certificate}}
\end{figure}

To show the impact of infeasibility certificate proposed in Section \ref {Infeasibility Certificate}, we conduct an experiment to compare the effects of using this certificate versus not using it, which is illustrated in Figure \ref{Sensitivity analysis for infeasibility certificate}. As shown in the figure, the total runtime is reduced by 5-15\% when using the infeasibility certificate. The average number of the reduction of the total runtime is equal to 11\%. This reduction is attributed to its ability to improve the efficiency of solving the Benders subproblem. In the experiment, the average and the maximum number of the reduction of GRM runtime are equal to 21\% and 28\% over these instances.

\subsubsection*{Customized initialization of gate patterns}
\begin{figure}[!htbp]
\centering
\includegraphics[width=0.75\textwidth]{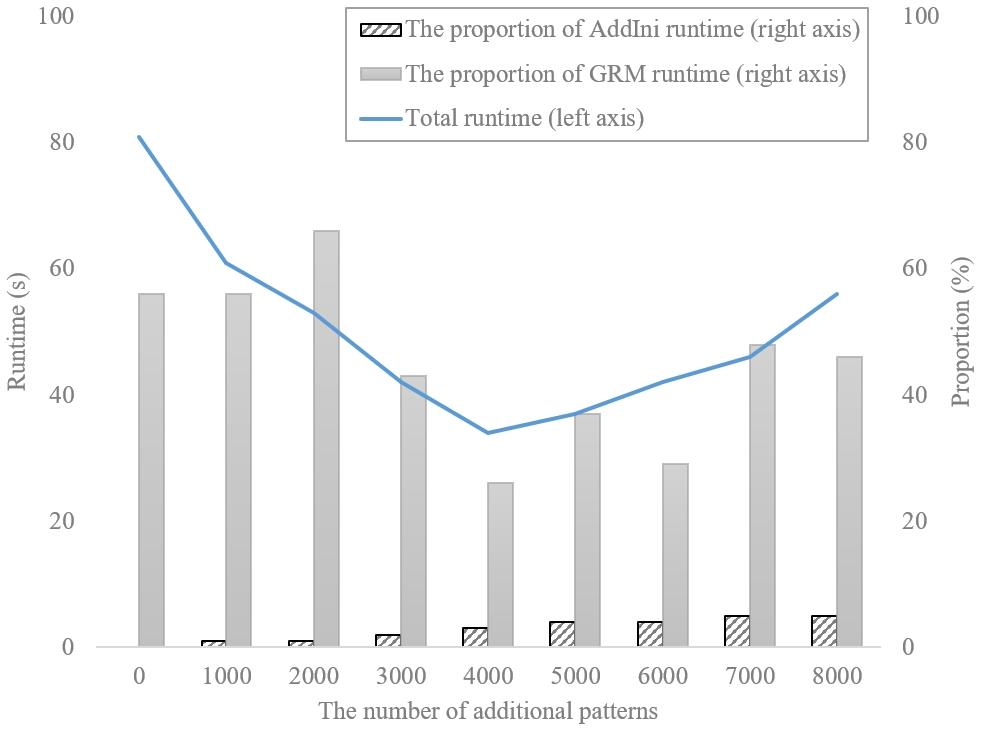}\\
\vspace{-0.15cm}
\hspace{-0.5cm}
{\footnotesize \textit{Note:} ``AddIni'' denotes additional initialization of gate patterns in solving RtBSP with CG}
\caption{The effect of additional gate pattern generation on computational time\label{Sensitivity analysis for additional patterns initialization}}
\end{figure}

There is a trade-off between the number of additional customized gate patterns (Section \ref{Section Benders subproblem (BSP)}) and the total runtime of the BCG method. We use the result from Instance \#2 to demonstrate this trade-off relationship. Other instances in the computational experiment have a similar trend of results. As shown in Figure \ref{Sensitivity analysis for additional patterns initialization}, with the increase of additional gate patterns initialized in solving the RxBSP, the extra runtime required for generating additional patterns becomes longer. Still, the runtime of solving the GRM becomes shorter. In addition, we find that the total runtime of the BCG method has a convex parabolic relationship with the number of additional patterns. When the number of additional patterns increases, the total runtime of the BCG method decreases at first and then gradually increases, with the minimum observed when adding an additional 4,000 gate patterns. The total runtime reduction is about 50\% by generating an additional 4,000 gate patterns compared to adding no extra gate patterns. By using an appropriate number of adding additional gate patterns, we can enhance the efficiency of the BCG method.\par

\subsubsection*{Optimality cuts}
Table \ref{Performance of the combination of optimality cuts} shows the effect of different combinations of the optimality cuts mentioned in Section \ref{Section Benders subproblem (BSP)} on Instance \#2. It is noted that Benders feasibility cuts and no-good cuts are all used in these combinations. From this table, we have the following observations. First, when we only employ the Benders optimality cuts (LPBOC), the optimality gap cannot reach below 5\% within 30 minutes. However, the optimality gaps corresponding to the combinations ``LPBOC+LLC'' and ``LPBOC+GC'' reach to 1.81\% and 3.85\% within 3 minutes, respectively. Second, according to the computational results, LLC seems effective for our problem, and the combination of LLC with other cuts can further reduce the optimality gap. For example, combining LLC with GC can reduce the optimality gap from 4.66\% with only GC to 0.84\% with ``LLC+GC''. Therefore, the results here further demonstrate the combinations of different optimality cuts can improve the solution quality of SAGRM and enhance the efficiency of the proposed BCG method.\par

\begin{table}[!htbp]
  \centering
  \renewcommand{\baselinestretch}{1.2}
  \caption{The performance of various combinations of optimality cuts\label{Performance of the combination of optimality cuts}}
  \small
    \begin{tabular}{lccc}
    \toprule
    Combination of cuts & Total recovery cost (unit) & Optimality gap (\%) & CPU time (s) \\
    \midrule
    LPBOC\footnotemark[1] & 199   & 35.14\% & 1800+\footnotemark[4] \\
    LLC\footnotemark[2]   & 191   & 4.86\% & 67 \\
    GC\footnotemark[3]    & 195   & 4.66\% & 396 \\
    LPBOC + LLC & 191   & 1.81\% & 51 \\
    LPBOC + GC & 195   & 3.85\% & 135 \\
    LLC + GC & 191   & 0.84\% & 48 \\
    LPBOC + LLC + GC & 191   & 0.70\% & 34 \\
    \bottomrule
    \end{tabular}\\
     \raggedright
     \vspace{0.1cm}
     \hspace{0.7cm}
     \footnotemark[1]{\footnotesize LPBOC: Benders optimality cut.}
     \hspace{0.1cm}
     \footnotemark[2]{\footnotesize LLC: Laporte \& Louveaux cut.}
     \hspace{0.1cm}
     \footnotemark[3]{\footnotesize GC: global cut.}\\
     \vspace{-0.1cm}
     \hspace{0.7cm}
     \footnotemark[4]1800+: out of the maximum allowable runtime, 1800 seconds.
\end{table}

\subsection{The effect of the length of delay interval and the number of gate types}\label{Section The length of delay interval and gate types}
The scale of the recovery problem is closely related to the length of delay interval and the number of gate types. On the one hand, a smaller interval will generate more flight delay options, but it will increase computational time. Hence, airlines must select an appropriate length of delay interval according to operational requirements. We examine this effect by testing different lengths of delay interval on Instance \#2. On the other hand, in order to highlight the complexity of available gate types, sensitivity analysis is done to explore the effect of the number of gate types as follows. We compare a scenario with one gate type and the other with two gate types (wide and narrow-body aircraft). In the scenario with one single gate type, every aircraft can park at any gate without the type restriction (Section \ref{Section Definitions} (6-a)). In the scenario with two gate types, aircraft should park at the gates whose types are consistent with the fleet type. Results are given in Table \ref{Results of different delay intervals and the number of gate types}.\par

\begin{table}[!htbp]
  \centering
  \renewcommand{\baselinestretch}{1}
  \caption{Scenario study results with different delay intervals and the number of gate types\label{Results of different delay intervals and the number of gate types}}
  \small
  {\resizebox{0.9\linewidth}{!}{
    \begin{tabular}{l|cccccccccc}
    \toprule
    Delay interval (mins) & \multicolumn{2}{c}{5} & \multicolumn{2}{c}{10} & \multicolumn{2}{c}{15} & \multicolumn{2}{c}{30} & \multicolumn{2}{c}{60} \\
    No. of flight copies & \multicolumn{2}{c}{1301} & \multicolumn{2}{c}{677} & \multicolumn{2}{c}{469} & \multicolumn{2}{c}{261} & \multicolumn{2}{c}{157} \\
    No. of gate types & 1     & 2     & 1     & 2     & 1     & 2     & 1     & 2     & 1     & 2 \\
    Total recovery cost (unit)  & 179   & 191   & 370   & 378   & 574   & 582   & 944   & 960   & 1480  & 1492  \\
    CPU time (s) & 25    & 34    & 19    & 27    & 15    & 21    & 11    & 18    & 8     & 12  \\
    \bottomrule
    \end{tabular}}}
\end{table}

As shown in Table \ref{Results of different delay intervals and the number of gate types}, the total recovery cost grows rapidly when the length of the delay interval increases. A longer interval for potential flight delay options could lead to fewer flight delay options (fewer flight copies) for each flight. In other words, the solution space of the recovery problem greatly narrows when the length becomes larger; hence, the total recovery cost increases. It is observed that the total recovery cost corresponding to a 60-minute delay interval is approximately seven times higher than the scenario with a 5-minute delay interval. On the other hand, the computational time decreases as expected when the length of the delay interval is longer.\par

Considering the effect of the number of gate types according to Table \ref{Results of different delay intervals and the number of gate types}, we find that the recovery cost of two gate types is always higher than that of one gate type. This is because gate reassignment must consider aircraft types when there are two types of gates. However, in the case that only one gate type is available, every aircraft can park at any gate without the gate type restriction. This makes the feasible set of two gate types a subset of the feasible set of one gate type. Thus, the total recovery cost associated with two gate types will be higher.\par

\section{Conclusion and future research}
\noindent
This paper considers the airline recovery problem, integrating flight rescheduling, aircraft rerouting, and gate reassignment to capture the interdependence among these recovery phases. A path-based model is formulated to describe the integrated recovery problem, where gate patterns are introduced to establish the relationship between flights and gates. Given the multi-stage characteristics of the integrated model, we propose a BCG method that consists of the Bender master problem for flight rescheduling and aircraft rerouting and the Benders subproblem for gate reassignment. To tackle a large number of decision variables, CG is used to generate aircraft routes and gate patterns iteratively. Based on the structure of our model, two acceleration techniques (a separation technique and an infeasibility certificate) are provided to enhance the BCG method. The computational experiments show that the proposed method can effectively solve the integrated recovery problem.\par

This work motivates us to consider future research in several directions. First, other recovery options (e.g., cruise speed control) can be considered to save the recovery cost further. Second, a customized parallel computing framework can be explored to promote the efficiency of the BCG method. Third, we would consider the uncertainty in airline disruption management to improve the practicality of our integrated recovery approach.\par

\setlength{\bibsep}{0.0em}

\newpage
\appendix
\section{The constraint programming models}\label{appendix CP models}
\subsection{The constraint programming model of aircraft routes}\label{The Constraint Programming Model of Aircraft Routes}
The constraint programming model for generating aircraft routes for each aircraft $r \in R$ is presented as follows. In this model, the decision variable $p=(p_{1}, p_{2}, ... , p_{\eta}) \in \{lb_r, lb_r + 1, ..., ub_r\}^{\eta}$ is a vector of the indexes of flight copies that can be performed by aircraft $r$, where $\eta$ denotes the maximum allowable number of legs in routes. Table \ref{Parameters used in the CP model1} provides the involved parameters. Then, we present the constraints to find possible aircraft routes.\par

\begin{table}[!htbp]
  \centering
  \renewcommand{\baselinestretch}{1}
  \caption{Parameters used in the CP model of aircraft routes}\label{Parameters used in the CP model1}
  \small
   \setlength{\tabcolsep}{4mm}{
    \begin{tabular}{ll}
    \toprule
    $DepAp_{j_v}/ArrAp_{j_v}$ & the departure/arrival airport of flight copy $j_v$\\
    $DepAp/ArrAp$ & the array of~$DepAp_{j_v}/ArrAp_{j_v}$\\
    $StartAP_{r}/EndAP_{r}$ & the start/end airport of aircraft $r$\\
    $DepT_{j_v}/ArrT_{j_v}/FC\_DelayT_{j_v}$ & the departure/arrival/delay time of flight copy $j_v$\\
    $DepT/ArrT/FC\_DelayT$ & the array of~$DepT_{j_v}/ArrT_{j_v}/DelayT_{j_v}$\\
    $TT_{r}$ & the minimum turn time of aircraft $r$\\
    $MT_{r}$ & the index of the maintenance task in $FC$ for aircraft $r$\\
    $UnplannedF_{rj}$ &1, if flight $j$ is unplanned for aircraft $r$;~0, otherwise\\
    $UnplannedF_{r}$ & the array of $UnplannedF_{rj}$ for aircraft $r$\\
    $R\_Dual_{r}$ & the dual value of aircraft $r$\\
    $FC\_Dual_{j_v}$ & the dual value related to flight copy $j_v$\\
    $FC\_Dual$ & the array of~$FC\_Dual_{j_v}$\\
    $Pair\_Dual_{j_{1v},j_{2v}}$ & the dual value related to the pair of flight copies  $(j_{1v},j_{2v})$\\
    $Pair\_Dual$ & the array of~$Pair\_Dual_{j_{1v},j_{2v}}$\\
    \bottomrule
    \end{tabular}
    }
\end{table}

\textbf{Constraints:}
\vspace{-0.5cm}
\begin{eqnarray}
&presenceOf(p_{k+1}) \leq presenceOf(p_{k}),~~\forall k\in \{1, ..., \eta-1\},\label{CP1 con1}\\
&presenceOf(p_{k+1}) == (p_{k+1}!=p_{k}),~~\forall k\in \{1, ..., \eta-1\},\label{CP1 con2}\\
&element(DepAP, p_1) == StartAP_{r},\label{CP1 con3}\\
&element(ArrAP, p_{\eta}) == EndAP_{r},\label{CP1 con4}\\
&\left(element(ArrAP, p_{k}) == element(DepAP, p_{k+1})\right)\notag\\
&\hspace*{5cm}\bigvee\left(presenceOf(p_{k+1}) == 0\right),~~\forall k\in \{1, ..., \eta-1\},\label{CP1 con5}\\
&\left(element(ArrT, p_{k}) + {TT_r} \leq element(DepT, p_{k+1})\right)\notag\\
&\hspace*{5cm}\bigvee\left(presenceOf(p_{k+1}) == 0\right),~~\forall k\in \{1, ..., \eta-1\},\label{CP1 con6}\\
&count(p, MT_{r}) \geq 1,~\textrm{if aircraft}~r~\textrm{has a maintenance task,}\label{CP1 con7}\\
&count(p, MT_{\bar{r}}) == 0,~~\forall \textrm{aircraft}~\bar{r}\in R\setminus\{r\}~\textrm{has a maintenance task,}\label{CP1 con8}\\
&\textrm{DelayT} = \sum\limits_{k=1}^{\eta} element(FC\_DelayT, p_{k})\cdot presenceOf(p_{k}), \label{CP1 con9}\\
&\textrm{SwapN} = \sum\limits_{k=1}^{\eta} element(UnplannedF_{r}, p_{k})\cdot presenceOf(p_{k}), \label{CP1 con10}\\
&c_p^r = c_{\textrm{Delay}} \cdot \textrm{DelayT} + c_{\textrm{Swap}} \cdot \textrm{SwapN}, \label{CP1 con11}\\
\nonumber
&\bar{c}_p^r = c_p^r + R\_Dual_r + \sum\limits_{k=1}^{\eta} element(FC\_Dual, p_{k})\cdot presenceOf(p_{k})\\
&~~~~+ \sum\limits_{k=1}^{\eta-1} element(Pair\_Dual, (p_{k},p_{k+1}))\cdot presenceOf(p_{k+1}). \label{CP1 con12}
\end{eqnarray}\par

The function $presenceOf(p_{k})$ indicates whether the element $p_{k}$ is present in the vector $p$. Suppose the value of $p_{k}$ is equal to $\alpha$, then the function $element(X, p_{k})$ returns the $\alpha$th element of the array $X$. Accordingly, suppose the values of $p_{k}$ and $p_{k+1}$ are equal to $\alpha$ and $\beta$ respectively, then the function $element(Y, (p_{k}, p_{k+1}))$ returns the ($\alpha$th, $\beta$th) element of the two-dimensional array $Y$. The function $count(p, \alpha)$ returns the number of the elements that are equal to the value $\alpha$ in the vector $p$.\par

Constraints \eqref{CP1 con1} restrict that $p_{k}$ must be present if $p_{k+1}$ is present. Constraints \eqref{CP1 con2} ensure that the value of $p_{k+1}$ is different from that of $p_{k}$ if $p_{k+1}$ is present; otherwise, the value of $p_{k+1}$ is the same as that of $p_{k}$. Constraints \eqref{CP1 con3} and \eqref{CP1 con4} restrict the start and end airports for aircraft $r$. Constraints \eqref{CP1 con5} and \eqref{CP1 con6} make sure that the space match and the minimum turn time are satisfied for each present flight connection in the vector $p$. Constraints \eqref{CP1 con7} and \eqref{CP1 con8} guarantee the maintenance requirement for aircraft $r$. Equations \eqref{CP1 con9}$-$\eqref{CP1 con11} calculate the cost of the aircraft route corresponding to the vector $p$. Equation \eqref{CP1 con12} calculates the reduced cost of the aircraft route corresponding to the vector $p$, where $R\_Dual_{r}, FC\_Dual$ and $Pair\_Dual$ can be obtained according to Equation \eqref{routeReducedCost}.\par

\subsection{The constraint programming model of gate patterns}\label{The Constraint Programming Model of Gate Patterns}
The constraint programming model of the gate patterns for the gates related to airport $a \in A$ and gate type $t \in T_a$ is presented as follows. In this model, the decision variable $p=(p_{1}, p_{2}, ... , p_{\eta}) \in \{lb_{a,t}, lb_{a,t}+1, ..., ub_{a,t}\}^{\eta}$ is a vector of the indexes of flight copies that can be serviced by the gates (i.e., satisfying the match constraint), where $\eta$ denotes the maximum allowable number of legs in patterns. The involved parameters have been given in Table \ref{Parameters used in the CP model1}. Now, we present the constraints to find possible gate patterns.\par

\vspace{0.2cm}
\textbf{Constraints:}
\vspace{-0.5cm}
\begin{eqnarray}
&presenceOf(p_{k+1}) \leq presenceOf(p_{k}),~~\forall k\in \{1, ..., \eta-1\},\label{CP2 con1}\\
&presenceOf(p_{k+1}) == (p_{k+1}!=p_{k}),~~\forall k\in \{1, ..., \eta-1\},\label{CP2 con2}\\
\nonumber
&T_k + \textrm{BufferT} \leq T_{k+1},~~\forall k\in \{1, ..., \eta-1\},~\textrm{where}\\
\nonumber
&T_k = element(DepT, p_k) \cdot \left[element(DepAP, p_k) == a\right] \cdot [element(DepAP, p_k)\\
\nonumber
&\hspace*{1.8cm}!= element(ArrAP, p_k)] + element(ArrT, p_k) \cdot \left[element(ArrAP, p_k) == a\right],\\
\nonumber
&\hspace*{-1cm}T_{k+1} = element(DepT, p_{k+1}) \cdot \left[element(DepAP, p_{k+1}) == a\right]  +  element(ArrT, p_{k+1})\cdot\\
&\hspace*{-0.2cm}  [element(ArrAP, p_{k+1}) == a] \cdot [element(DepAP, p_{k+1}) != element(ArrAP, p_{k+1})]\label{CP2 con3}\\
\nonumber
&\bar{c}_p^{a,t} =  c_{a,t} + G\_Dual_{a,t} - \sum\limits_{k=1}^{\eta} element(FC\_Dual, p_{k})\cdot presenceOf(p_{k})\\
&~~ - \sum\limits_{k=1}^{\eta-1} element(Pair\_Dual, (p_{k},p_{k+1}))\cdot presenceOf(p_{k+1}). \label{CP2 con4}
\end{eqnarray}\par

Constraints \eqref{CP2 con1}$-$\eqref{CP2 con2} are similar to Constraints \eqref{CP1 con1}$-$\eqref{CP1 con2}. Constraints \eqref{CP2 con3} ensure the buffer time constraint between two consecutive activities, where $\textrm{BufferT}$ denotes the value of the buffer time. The terms $[element(DepAP, p_k) != element(ArrAP, p_k)]$ and $[element(DepAP, p_{k+1}) != element(ArrAP, p_{k+1})]$ are introduced to avoid repetitively calculating the start/end time of maintenance activities. Equation \eqref{CP2 con4} calculates the reduced cost of the gate pattern corresponding to the vector $p$, where $G\_Dual_{a,t}, FC\_Dual$ and $Pair\_Dual$ can be obtained according to Equation \eqref{patternReducedCost}.\par


\addcontentsline{toc}{section}{References}
\begin{thebibliography}{60}
\bibitem[Ahmed et al.(2018)]{ahmed2018robust} Ahmed, M.B., Mansour, F.Z., Haouari, M., 2018. Robust integrated maintenance aircraft routing and crew pairing. J. Air Transport Manag. 73, 15-31.
\bibitem[Apt(2003)]{apt2003} Apt, K., 2003. Principles of Constraint Programming. Cambridge University Press, Cambridge.
\bibitem[Ar{\i}kan et al.(2016)]{arikan2016integrated} Ar{\i}kan, U., G{\"u}rel, S., Akt{\"u}rk, M.S., 2016. Integrated aircraft and passenger recovery with cruise time controllability. Annu. Oper. Res. 236(2), 295-317.
\bibitem[Ar{\i}kan et al.(2017)]{arikan2017flight} Ar{\i}kan, U., G{\"u}rel, S., Akt{\"u}rk, M.S., 2017. Flight network-based approach for integrated airline recovery with cruise speed control. Transp. Sci. 51(4), 1259-1287.
\bibitem[Barnhart et al.(1998)]{barnhart1998flight} Barnhart, C., Boland, N.L., Clarke, L.W., Johnson, E.L., Nemhauser, G.L., Shenoi, R.G., 1998. Flight string models for aircraft fleeting and routing. Transp. Sci. 32(3), 208-220.
\bibitem[Benoist et al.(2002)]{benoist2002constraint} Benoist, T., Gaudin, E., Rottembourg, B., 2002. Constraint programming contribution to Benders decomposition: A case study. Principles and Practice of Constraint Programming-CP 2002, Springer, Berlin, 603-617.
\bibitem[Boyd and Vandenberghe(2004)]{boyd2004convex} Boyd, S., Vandenberghe, L., 2004. Convex Optimization, Cambridge University Press, UK.
\bibitem[Bratu and Barnhart(2006)]{bratu2006flight} Bratu, S., Barnhart, C., 2006. Flight operations recovery: New approaches considering passenger recovery. J. Sched. 9, 279-298.
\bibitem[Brueckner et al.(2021)]{BRUECKNER2021102333} Brueckner, J.K., Czerny, A.I., Gaggero, A.A., 2021. Airline mitigation of propagated delays via schedule buffers: Theory and empirics. Transp. Res. Part E: Logist. Transp. Rev. 150, 102333.
\bibitem[Cadarso and Vaze(2023)]{cadarso2022} Cadarso, L., Vaze, V., 2023. Passenger-centric integrated airline schedule and aircraft recovery. Transp. Sci. 57(3), 813-837.
\bibitem[Castaing et al.(2016)]{castaing2016reducing} Castaing, J., Mukherjee, I., Cohn, A., Hurwitz, L., Nguyen, A., M{\"u}ller, J.J., 2016. Reducing airport gate blockage in passenger aviation: Models and analysis. Comput. Oper. Res. 65, 189-199.
\bibitem[Clausen et al.(2010)]{clausen2010disruption} Clausen, J., Larsen, A., Larsen, J., Rezanova, N.J., 2010. Disruption management in the airline industry$-$Concepts, models and methods. Comput. Oper. Res. 37(5), 809-821.
\bibitem[Cordeau et al.(2001)]{cordeau2001benders} Cordeau, J.F., Stojkovi{\'c}, G., Soumis, F., Desrosiers, J., 2001. Benders decomposition for simultaneous aircraft routing and crew scheduling. Transp. Sci. 35(4), 375-388.
\bibitem[Diepen et al.(2012)]{diepen2012gate} Diepen, G., van den Akker, J.M., Hoogeveen, J.A., Smeltink, J.W., 2012. Finding a robust assignment of flights to gates at Amsterdam Airport Schiphol. J. Sched. 15, 703-715.
\bibitem[Ding et al.(2023)]{ding2023towards} Ding, Y., Wandelt, S., Wu, G., Xu, Y., Sun, X., 2023. Towards efficient airline disruption recovery with reinforcement learning. Transp. Res. Part E: Logist. Transp. Rev. 179, 103295.
\bibitem[Dorndorf et al.(2017)]{dorndorf2017reducing} Dorndorf, U., Jaehn, F., Pesch, E., 2017. Flight gate assignment and recovery strategies with stochastic arrival and departure times. OR Spect. 39, 65-93.
\bibitem[Dunbar et al.(2012)]{dunbar2012robust} Dunbar, M., Froyland, G., Wu, C.L., 2012. Robust airline schedule planning: Minimizing propagated delay in an integrated routing and crewing framework. Transp. Sci. 46(2), 204-216.
\bibitem[Fakhri et al.(2017)]{fakhri2017benders} Fakhri, A., Ghatee, M., Fragkogios, A., Saharidis, G.K.D., 2017. Benders decomposition with integer subproblem. Expert Syst. Appl. 89, 20-30.
\bibitem[Flyertalk(2019)]{flyertalk2019} Flyertalk, 2019. The ``reverse delay" (gate not available on arrival) at ORD \& elsewhere, accessed February 1, 2025, https://www.flyertalk.com/forum/united-airlines-mileageplus/1949747-reverse-delay-gate-not-available-arrival-ord-elsewhere.html.
\bibitem[Geske et al.(2024)]{geske2024} Geske, A.M, Herold, D.M., Kummer, S., 2024. Integrating AI support into a framework for collaborative decision-making (CDM) for airline disruption management. J. Air Transp. Res. Soc. 3, 100026.
\bibitem[Hassan et al.(2021)]{hassan2021airline} Hassan, L.K., Santos, B.F., Vink, J., 2021. Airline disruption management: A literature review and practical challenges. Comput. Oper. Res. 127, 105137.
\bibitem[Hu et al.(2016)]{hu2016integrated} Hu, Y., Song, Y., Zhao, K., Xu, B., 2016. Integrated recovery of aircraft and passengers after airline operation disruption based on a GRASP algorithm. Transp. Res. Part E: Logist. Transp. Rev. 87, 97-112.
\bibitem[Hu et al.(2021)]{hu2021integrated} Hu, Y., Zhang, P., Fan, B., Zhang, S., Song, J., 2021. Integrated recovery of aircraft and passengers after airline operation disruption based on a GRASP algorithm. Comput. Ind. Eng. 161, 107664.
\bibitem[Jiang et al.(2025)]{JIANG2025104243} Jiang, J., Zhang, S., Tang, Y., Guo, Y., Wu, C.L, 2025. ADMM-based augmented Lagrangian methods for robust aircraft recovery problem considering connection time, resource capacity and maintenance flexibility. Transp. Res. Part E: Logist. Transp. Rev. 201, 104243.
\bibitem[Kontoyiannakis et al.(2009)]{kontoyiannakis2009simulation} Kontoyiannakis, K., Serrano, E., Tse, K., Lapp, M., Cohn, A., 2009. A simulation framework to evaluate airport gate allocation policies under extreme delay conditions. Proceedings of the 2009 Winter Simulation Conference (WSC), IEEE, 2332-2342.
\bibitem[Laporte and Louveaux(1993)]{laporte1993integer} Laporte, G., Louveaux, F.V., 1993. The integer L-shaped method for stochastic integer programs with complete recourse. Oper. Res. Lett. 13(3), 133-142.
\bibitem[Liang et al.(2018)]{liang2018column} Liang, Z., Xiao, F., Qian, X., Zhou, L., Jin, X., Lu, X., Karichery, S., 2018. A column generation-based heuristic for aircraft recovery problem with airport capacity constraints and maintenance flexibility. Transp. Res. Part B: Methodol. 113, 70-90.
\bibitem[Maharjan and Matis(2011)]{maharjan2011optimization} Maharjan, B., Matis, T.I., 2011. An optimization model for gate reassignment in response to flight delays. J. Air Transport Manag. 17(4), 256-261.
\bibitem[Maher(2015)]{maher2015novel} Maher, S.J., 2015. A novel passenger recovery approach for the integrated airline recovery problem. Comput. Oper. Res. 57, 123-137.
\bibitem[Maher(2016)]{maher2016solving} Maher, S.J., 2016. Solving the integrated airline recovery problem using column-and-row generation. Transp. Sci. 50(1), 216-239.
\bibitem[Papadakos(2009)]{papadakos2009integrated} Papadakos, N., 2009. Integrated airline scheduling. Comput. Oper. Res. 36(1), 176-195.
\bibitem[Petersen et al.(2012)]{petersen2012optimization} Petersen, J.D., S{\"o}lveling, G., Clarke, J.P., Johnson, E.L., Shebalov, S., 2012. An optimization approach to airline integrated recovery. Transp. Sci. 46(4), 482-500.
\bibitem[Poyraz and Azizo\v{g}lu(2024)]{poyraz2024} Poyraz, D.D., Azizo\v{g}lu, M., 2024. An airport gate reassignment problem with gate closures. J. Air Transport Manag. 115, 102529.
\bibitem[Pternea and Haghani(2019)]{pternea2019aircraft} Pternea, M., Haghani, A., 2019. An aircraft-to-gate reassignment framework for dealing with schedule disruptions. J. Air Transport Manag. 78, 116-132.
\bibitem[PYOK(2024)]{PYOK2024} PYOK, 2024. Man gets maximum prison sentence after assaulting a Southwest flight attendant when there wasn¡¯t a gate ready after arrival, accessed February 1, 2025, https://www.paddleyourownkanoo.com/2024/02/01/man-gets-maximum-prison-sentence-after-assaulting-a-southwest-flight-attendant-when-there-wasnt-a-gate-ready-after-arrival/.
\bibitem[Rosenberger et al.(2003)]{rosenberger2003rerouting} Rosenberger, J.M., Johnson, E.L., Nemhauser, G.L., 2003. Rerouting aircraft for airline recovery. Transp. Sci. 37(4), 408-421.
\bibitem[Rossi et al.(2006)]{rossi2006} Rossi, F., van Beek, P., Walsh, T., 2006. Handbook of Constraint Programming. Elsevier, Amsterdam.
\bibitem[Sinclair et al.(2014)]{sinclair2014} Sinclair, K., Cordeau, J.F., Laporte, G., 2014. Improvements to a large neighborhood search heuristic for an integrated aircraft and passenger recovery problem. Eur. J. Oper. Res. 233(1), 234-245.
\bibitem[Sinclair et al.(2016)]{sinclair2016} Sinclair, K., Cordeau, J.F., Laporte, G., 2016. A column generation post-optimization heuristic for the integrated aircraft and passenger recovery problem. Comput. Oper. Res. 65, 42-52.
\bibitem[Su et al.(2021)]{su2021airline} Su, Y., Xie, K., Wang, H., Liang, Z., Chaovalitwongse, W.A., Pardalos, P.M., 2021. Airline disruption management: A review of models and solution methods. Engineering 7(4), 435-447.
\bibitem[Tang and Wang(2013)]{tang2013airport} Tang, C.H., Wang, W.C., 2013. Airport gate assignments for airline-specific gates. J. Air Transport Manag. 30, 10-16.
\bibitem[Teodorovi{\'c} and Guberini{\'c}(1984)]{teodorovic1984optimal} Teodorovi{\'c}, D., Guberini{\'c}, S., 1984. Optimal dispatching strategy on an airline network after a schedule perturbation. Eur. J. Oper. Res. 15(2), 178-182.
\bibitem[Teodorovi{\'c} and Stojkovi{\'c}(1990)]{teodorovic1990model} Teodorovi{\'c}, D., Stojkovi{\'c}, G., 1990. Model for operational daily airline scheduling. Transp. Plan. Technol. 14(4), 273-285.
\bibitem[Teodorovi{\'c} and Stojkovi{\'c}(1995)]{teodorovic1995model} Teodorovi{\'c}, D., Stojkovi{\'c}, G., 1995. Model to reduce airline schedule disturbances. J. Transp. Eng. 121(4), 324-331.
\bibitem[Thengvall et al.(2001)]{thengvall2001multiple} Thengvall, B.G., Yu, G., Bard, J.F., 2001. Multiple fleet aircraft schedule recovery following hub closures. Transp. Res. Part A: Policy Pract. 35(4), 289-308.
\bibitem[Wang et al.(2019)]{wang2019} Wang, D., Wu, Y., Hu, J.Q, Liu, M., Yu, P., Zhang, C., Wu, Y., 2019. Flight schedule recovery: A simulation-based approach. Asia-Pac. J. Oper. Res. 36(6), 1940010.
\bibitem[Wang et al.(2025)]{wang2025} Wang, Q., Mao, J., Wen, X., Wallace, S.W., Deveci, M., 2025. Flight, aircraft, and crew integrated recovery policies for airlines--A deep reinforcement learning approach. Transport Policy. 160, 245-258.
\bibitem[WLOS(2024)]{WLOS2024} WLOS, 2024. Passengers left waiting on Asheville Regional Airport tarmac for about two hours, accessed February 1, 2025, https://wlos.com/news/local/passengers-left-waiting-on-asheville-regional-airport-tarmac-for-about-two-hours.
\bibitem[Wolsey(2020)]{wolsey2020integer} Wolsey, L.A., 2020. Integer Programming, 2nd ed., John Wiley \& Sons, New York.
\bibitem[Wu and Maher(2018)]{wu2018airlinecapacity} Wu, C.L., Maher, S.J., 2018. Airline capacity planning and management, in The Routledge Companion to Air Transport Management ed. Nigel Halpern and Anne Graham, Abingdon: Routledge, January 31, 2018, accessed February 1, 2025, Routledge Handbooks Online.
\bibitem[Wu and Law(2019)]{WU201962} Wu, C.L., Law, K., 2019. Modelling the delay propagation effects of multiple resource connections in an airline network using a Bayesian network model. Transp. Res. Part E: Logist. Transp. Rev. 122, 62-77.
\bibitem[Wu et al.(2025)]{wu2025} Wu, S., Liu, E., Cao, R., Bai, Q., 2025. Airline recovery problem under disruptions: A review. Comput. Oper. Res. 175, 106915.
\bibitem[Xiong and Hansen(2013)]{XIONG201364} Xiong, J., Hansen, M., 2013. Modelling airline flight cancellation decisions. Transp. Res. Part E: Logist. Transp. Rev. 56, 64-80.
\bibitem[Xu et al.(2023)]{xu2023distributionally} Xu, Y., Wandelt, S., Sun, X., 2023. A distributionally robust optimization approach for airline integrated recovery under in-flight pandemic transmission risks. Transp. Res. Part C: Emerg. Technol. 152, 104188.
\bibitem[Yan et al.(2016)]{yan2016tarmac} Yan, C., Vaze, V., Vanderboll, A., Barnhart, C., 2016. Tarmac delay policies: A passenger-centric analysis. J. Air Transport Manag. 30, 10-16.
\bibitem[Yan et al.(2009)]{yan2009airport} Yan, S., Chen, C.Y., Tang, C.H., 2009. Airport gate reassignment following temporary airport closures. Transportmetrica 5(1), 25-41.
\bibitem[Yan et al.(2011)]{yan2011airport} Yan, S., Tang, C.H., Hou, Y.Z., 2011. Airport gate reassignments considering deterministic and stochastic flight departure/arrival times. J. Adv. Transp. 45(4), 304-320.
\bibitem[Yan and Yang(1996)]{yan1996decision} Yan, S., Yang, D.H., 1996. A decision support framework for handling schedule perturbation. Transp. Res. Part B: Methodol. 30(6), 405-419.
\bibitem[Zang et al.(2024)]{zang2024} Zang, H., Zhu, J., Zhu, Q., Gao, Q., 2024. A proactive aircraft recovery approach based on airport spatiotemporal network supply and demand coordination. Comput. Oper. Res. 165, 106599.
\bibitem[Zhang and Klabjan(2017)]{zhang2017optimization} Zhang, D., Klabjan, D., 2017. Optimization for gate re-assignment. Transp. Res. Part B: Methodol. 95, 260-284.
\bibitem[Zhang et al.(2016)]{zhang2016math} Zhang, D., Yu, C., Desai, J., Lau, H.Y.K., 2016. A math-heuristic algorithm for the integrated air service recovery. Transp. Res. Part B: Methodol. 84, 211-236.
\bibitem[Zhong et al.(2024)]{zhong2024} Zhong, H., Lian, Z., Zhou, T., Niu, B., 2024. A time-varying competitive swarm optimizer for integrated flight recovery with multi-objective and priority considerations. Comput. Ind. Eng. 190, 110019.
\end{thebibliography}
\end{document}